\documentclass[12pt,a4paper]{amsart}

\usepackage[english]{babel}
\usepackage[utf8]{inputenc}

\makeatletter
\@namedef{subjclassname@2022}{\emph{2022} Math Subject Classification}
\makeatother

\theoremstyle{definition} 
\newtheorem{definicija}{Definition}[section]
\newtheorem{primer}[definicija]{Example}
\newtheorem{opomba}[definicija]{Remark}

\theoremstyle{plain} 
\newtheorem{lema}[definicija]{Lemma}
\newtheorem{izrek}[definicija]{Theorem}
\newtheorem{trditev}[definicija]{Proposition}
\newtheorem{posledica}[definicija]{Corollary}

\DeclareMathOperator{\sgn}{sgn}
\DeclareMathOperator{\supp}{supp}

\newcommand{\R}{\mathbb R}
\newcommand{\N}{\mathbb N}
\newcommand{\Z}{\mathbb Z}
\newcommand{\C}{\mathbb C}
\newcommand{\Q}{\mathbb Q}

\newcommand{\cweyl}{\mathcal{A}_1(\C)}
\newcommand{\rweyl}{\mathcal{A}_1(\R)}

\allowdisplaybreaks
\numberwithin{equation}{section}
\makeatletter
\@namedef{subjclassname@2020}{\textup{2020} Mathematics Subject Classification}
\makeatother
\title{Valuations and orderings on the real Weyl algebra}
\author[L. Vukšić]{Lara Vukšić}

\address{Lara Vukšić, Institute of Mathematics, Physics, and Mechanics, Ljubljana, Slovenia and Faculty of Mathematics and Physics, University of Ljubljana, Slovenia}

\email{lara.vuksic@fmf.uni-lj.si}

\thanks{ I would like to thank my advisor Igor Klep for his guidance and many helpful comments and suggestions.}


\date{\today}


\title{Valuations and orderings on the real Weyl algebra}
\begin{document}
\maketitle
\begin{abstract}
The first Weyl algebra $\mathcal{A}_1(k)$ over a field $k$ is the $k$-algebra with two generators $x, y$ subject to $[y, x] = 1$ and was first introduced during the development of quantum mechanics. In this article, we classify all valuations on the real Weyl algebra $\rweyl$ whose residue field is $\R$. We then use a noncommutative version of the Baer-Krull theorem to classify all orderings on $\rweyl$. As a byproduct of our studies, we settle two open problems in real algebraic geometry. First, we show that not all orderings on $\rweyl$ extend to an ordering on a larger ring $R[y ; \delta]$, where $R$ is the ring of Puiseux series, introduced by Marshall and Zhang in \cite{MA2}, and characterize the orderings that do have such an extension. Second, we show that for valuations on noncommutative division rings, Kaplansky's theorem that extensions by limits of pseudo-Cauchy sequences are immediate fails in general.
\end{abstract}

\tableofcontents

\section{Introduction}
Valuation theory was first developed for commutative fields in the context of number theory and was first defined by József Kürschák \cite{KU} in 1913. For modern treatments, we refer to the books of Engler and Prestel \cite{ENG} or Kuhlmann \cite{KUHL2}. Oscar Schilling wrote the first major work on valuations on (noncommutative) division rings in 1945 \cite{SCH}.  

A valuation on a division ring $D$ is a map $v : D \rightarrow \Gamma \cup \{ \infty \}$, where $\Gamma$ is an ordered group written additively and $\infty \not \in \Gamma, \infty > \gamma$ for each $\gamma \in \Gamma$, with the following properties:
\begin{enumerate}
\item $\forall x \in D: v(x) = \infty \Leftrightarrow x = 0$,
\item $\forall x, y \in D : v( xy ) = v(x) + v(y)$,
\item $\forall x, y \in D : v(x + y) \ge \min\{ v(x), v(y)  \}$.
\end{enumerate}
It follows that $v$ is a homomorphism from $D^*$ to $\Gamma$.
The set $\mathcal{O}_v := \{ x \in D \mid v(x) \ge 0  \}$ is called the valuation ring associated to $v$, and $\mathcal{M}_v :=  \{ x \in M \mid v(x) > 0   \}$ is its maximal ideal. The division ring $\overline{D} := \mathcal{O}_v / \mathcal{M}_v$ is called the associated residue division ring.  Since $v$ is a group homomorphism, the subgroup $\mathcal{O}_v^*$ is normal in $D^*$. 
Several alternative approaches to noncommutative valuations, where $v$ does not define a group homomorphism, were introduced and studied recently by 
Nicolai Ivanovich Dubrovin in \cite{DU1} and \cite{DU2} (see also \cite{MARU} for a more thorough treatment), and by Jean-Pierre Tignol and Adrien Wadsworth in \cite{TIGN}.

Suppose $F$ is a field. Then all valuations on the field or rational functions $F(x)$ with residue field $F$ are well-known, namely, the $p$-adic valuations for irreducible polynomials $p(x) \in F[x]$, and the $v_{\deg}$ valuation, defined by $$v_{\deg}(\frac{p}{q} ) := \deg( q) - \deg (p).$$ The description of all valuations on the field of rational functions in several variables with residue field equal to the base field is much more involved. There are many descriptions of constructions of such valuations in the literature. Among famous examples of such descriptions are the one given by Saunders MacLane in \cite{MAC} and the one given by Franz-Viktor Kuhlmann in \cite{KUHL}. 

As valuations on Ore extensions uniquely extend to their quotient division ring, the description of all valuations on Ore division rings is equivalent to the description of all valuations on corresponding quotient division rings. The description of all valuations on noncommutative Ore extensions $R[x; \sigma,  \delta  ]$ where $R$ is a domain, 
$\sigma: R \rightarrow R$ is a ring homomorphism and $\delta: R \rightarrow R$ a $\sigma$-derivation is even more complex than in the commutative case. 
Additional difficulties arise from the fact that $[f, g] = 0 $ does not hold for all $f, g \in R[x; \sigma,  \delta  ]$. 
Granja, Martínez, and Rodríguez have shown in \cite{GRA} that the set of all real valuations extending to the skew polynomial ring has the structure of a parameterized complete non-metric tree.
Further recent progress on valuations on Ore extensions is given by Onay in \cite{GON} and Rohwer in his PhD thesis \cite{ROH}.

\subsection{Results}
Our main goal is to classify all orderings and real valuations on the real Weyl algebra $\rweyl$ or, equivalently, its quotient division ring $\mathcal{D}_1(\R)$. The Weyl algebra is the noncommutative algebra generated by two elements $x, y$ satisfying $[y, x] = 1$. Hence its elements are all of the form 
$$ \sum_{i,j} \alpha_{i,j}x^iy^j, \  \alpha_{i,j} \in \R.
$$
Because of this, our approach to constructing valuations on $\rweyl$ is inspired by classical constructions  of valuations on commutative rational functions in two unknowns mentioned above. However, the relation $[y, x] = 1$ gives rise to additional constraints and many fewer valuations than in the commutative case.

As we will show, the valuations on $\rweyl$ we are interested in all satisfy $v[a, b] > v(ab)$ for all nonzero $a, b$. We call such valuations  \emph{strongly abelian}. They have an abelian value group and commutative residue field. 
In Section \ref{properties} we give some properties of strongly abelian valuations. We show that if a valuation $v$ on a division ring $D$ satisfies $\overline{D} = \overline{Z(D)}$ and the value group is of rational rank one, then $v$ is strongly abelian. Under additional constraints on the residue field and the value group we extend this statement to valuations of higher rational rank.

In Section \ref{rweyl} we give a characterization of all valuations $v$ on the real Weyl algebra $\rweyl$ with residue field $\R$ in the spirit of MacLane \cite{MAC}. 
The construction is inspired by the outline given by Shtipelman in \cite{SHT} for  valuations on the complex Weyl algebra $\cweyl$. 
We also explicitly describe the associated value groups and show that they are all isomorphic to subgroups of $\Q$ or $\Q \times \Z$.

In their attempt to describe all orderings on $\rweyl$ in \cite{MA2}, Murray Marshall and Yufei Zhang introduced the Ore extensions $R[y ; \delta]$ and $\tilde{R}[y ; \delta]$, with $$ R := \{ \sum_{k \ge m} a_k x^{- \frac{k}{n} }  \mid a_k \in \R, m \in \Z, n \in \N  \},$$
 $$ \tilde{R} := \{ \sum_{q \in A} a_q x^{- q }  \mid A \subset \Q \text{ is well-ordered} \}$$  and  $\delta(p(x)) = p'(x)$. 
As is often done in real algebraic geometry, all orderings are described by classifying all real valuations via the Baer-Krull theorem.
Marshall and Zhang described \emph{almost} all valuations $v$ on $R[y ; \delta]$ with residue field $\R$; in one case, they did not prove that $v$ is a valuation. In Section \ref{middle}, we complete their characterization. Marshall and Zhang also conjectured that all valuations on $\rweyl$ with residue field $\R$ extend to a valuation on $R[y ; \delta]$ with the same residue field.
We refute their conjecture in Section \ref{middle}. Further, we combine our classification of valuations on $\rweyl$ with
Marshall and Zhang's description of valuations on $R[y ; \delta]$ to characterize the valuations on $\rweyl$ with residue field $\R$ that extend to a valuation $R[y , \delta]$ with the same residue field. All such extensions are again strongly abelian. 

In Section \ref{grand}, we show that all valuations on $R[y ; \delta]$ with residue field $\R$ uniquely extend to a strongly abelian valuation on $\tilde{R}[y ; \delta]$ with the same residue field. We also show that the value group of such an extension is not of rational rank one. 
As a byproduct of our investigations, we show that Kaplansky's theorem that all extensions by limits of pseudo-Cauchy sequences are immediate (in particular, they do not change the rational rank of the value group) fails for noncommutative division rings.

As Marshall and Zhang observe in \cite{MA1}, all strongly abelian valuations $v$ on a division ring $D$ with a formally real residue field are compatible with an order on $D$. In Section \ref{orders}, we describe all $v$-compatible orders on $\rweyl$ for every valuation $v$ on $\rweyl$ constructed in Section \ref{rweyl} using a noncommutative version of the Baer-Krull theorem as given in \cite{CON} (see also
\cite{TSC}, \cite{CIM}, \cite{CRA} and \cite{KLE} for modern treatments and extensions). We also characterize the $v$-compatible orders on $\rweyl$ that extend to an order on $R[y ; \delta]$ compatible with $v$'s extension to $R[y ; \delta]$.

\section{Strongly abelian valuations}\label{properties}

We present some properties of valuations on noncommutative division rings which we will use later to describe order-compatible valuations on the real Weyl algebra $\rweyl$ and some of its ring extensions. First, we define a property of valuations on division rings.

\begin{definicija}
Suppose $v$ is a valuation on a division ring $D$. We say $v$ is \emph{strongly abelian} if $ v[a, b] > v(ab)$ holds for all nonzero $a, b \in D$.
\end{definicija}

Any valuation on a field is strongly abelian. In this section, we describe a sufficient condition for a valuation $v$ to be strongly abelian.
This property will be important for us for two reasons. Firstly, it is obvious than if a valuation $v$ on a division ring $D$ is strongly abelian, then the associated value group and residue division ring are commutative. Secondly, we are particularly interested in order-compatible valuations on $\rweyl$; minimal such have residue field $\R$, as it was shown in \cite{MA2}. It follows from Theorem 2.5 of \cite{MA1} that a strongly abelian valuation $v$ on a division ring $D$ with a formally real residue field is compatible with an order on $D$ by the noncommutative version of the Baer-Krull theorem as given in \cite{TSC}.

\begin{trditev} \label{ratdepelt}
Let $v$ be a valuation on a division ring $D$ such that $\overline{D} = \overline{Z(D)}$. Let $a, b \in D^*$ be such that $v(a)$ and
$v(b)$ are rationally dependent. Then $v[a, b] > v(ab)$.
\end{trditev}

\proof 
Since $v(a)$ and $v(b)$ are rationally dependent, $v(ab) = v(ba) \le v[a, b ]$. Suppose $v[a, b] = v(ab) $ and let $\beta := \overline{aba^{-1}b^{-1}} \in \overline{D}$. We have

\begin{align*} 
\beta = \overline{aba^{-1}b^{-1}} = \overline{( [a,b] + ba   )a^{-1}b^{-1}} = \overline{[a,b]a^{-1}b^{-1} + 1} \neq 1.
\end{align*}

Let $v(b) = - \frac{\ell}{k} v(a)$ for $\ell, k \in \Z$, $\ell$ and $k$ coprime. It follows that $v[a, b] = \frac{k - \ell}{k}v(a).$ Define $\gamma := \overline{ (ba  )^k a^{ \ell - k}   } \in \overline{D}$. Let $\beta', \gamma' \in Z(D)$ be such that $v(\beta' ) = v(\gamma') = 0$, $\overline{\beta'} = \beta$  and $\overline{\gamma'} = \gamma$. Then on one hand,

\begin{align*}
v(  a(ba)^k - \gamma' a^{k - \ell + 1}    ) = v( a ( (ba)^k a^{ \ell - k} - \gamma'  )a^{k- \ell}  ) > (k -  \ell + 1)v(a),
\end{align*}
and on the other, 
\begin{align*}
v(  a(ba)^k - \gamma' a^{k - \ell + 1}    ) &= v( (ab)^k a - \gamma' a^{k - \ell + 1}   ) = v( ((ab)^ka^{\ell - k}  - \gamma'  )a^{k - \ell + 1}  ) \\
&= (k - \ell + 1)v(a)
\end{align*}
since 
$$ \overline{ (ab)^k a^{  \ell - k} - \gamma'} = \overline{  (aba^{-1}b^{-1}ba    )^k a^{\ell-k}} - \gamma     =  \beta^k\overline{ ( ba)^k a^{\ell - k}} - \gamma   \neq 0   . $$  
In the last equation, we used that $\overline{(aba^{-1}b^{-1}ba    )^ka^{\ell-k}} = \overline{ ( \beta'  ba   )^ka^{\ell-k} } = \overline{\beta'^k(ba)^ka^{\ell-k}}  =  \beta^k \overline{ (ba)^k  a^{\ell-k}}$ since $\overline{D} = \overline{Z(D)}$ as presumed.
This holds if $\beta^k \neq 1$. If $\beta^k = 1$, we repeat our calculations with roles of $a$ and $b$ interchanged. The new $\beta$ will now be the inverse value of the former and since $\gcd (\ell, k) = 1$, $\beta^{- \ell } \neq 1$. In either case, we get a contradiction from which we deduce $v[a, b] > v(ab)$. 
\endproof

\begin{opomba}  The condition  $\overline{D} = \overline{Z(D)}$ is fulfilled by every valuation on an algebra over a field that is isomorphic to the residue field. In particular, this holds for minimal order-compatible valuations on $\R$-algebras.
\end{opomba}

\begin{posledica} \label{ratvalgrp}
Let $v$ be a valuation on a division ring $D$ such that $\overline{D} = \overline{Z(D)}$. If the value group has rational rank one, then $v$ is strongly abelian.
\end{posledica}
 
\begin{lema}\label{integerpowers}
Let $v$ be a valuation on a division ring $D$ such that the value group is abelian and $\overline{D} = \overline{Z(D)}$. Then for all $x, y \in D \setminus \{0 \}$:
\begin{enumerate}
\item If $v[x, y] > v(xy) $, then $v[x^m, y] > v( x^m y)$ for all $m \in \Z$.
\item Suppose $v[x, y] = v(xy) $. Then $v[x^{-1}, y] = v(x^{-1}y) $ and for each $m \in \N$,  $v[x^m, y]  >  v( x^m y)$ holds if and only if
  ${\alpha} := {y^{-1}x^{-1}yx}$ satisfies $1 + \overline{\alpha} + \cdots + \overline{\alpha}^{m - 1} = 0$ in $\overline{D}.$
\end{enumerate}
\end{lema}

\proof
To prove (1), first observe
$$[x^{-1}, y] = x^{-1}y - y x^{-1} = x^{-1}(yx - xy)x^{-1},$$ so if $v[x,y] > v(x, y)$, then
 $$v[x^{-1}, y] = v[x, y] - 2v(x) > v(xy) - 2v(x) = v(x^{-1}y). $$
Suppose $m \in \N$. Then
\begin{align*}
[ x^m, y  ] &=   \sum_{\ell = 1}^m x^{m - \ell} [x, y] x^{\ell - 1} 
\end{align*}
 and since the value group is commutative, $v(x^{m - \ell} [x, y] x^{\ell - 1}) = (m-1)v(x) + v[x, y] > v(x^my)$ for each $1 \le \ell \le m$. Item (1) is thus proved.

To prove (2), suppose $v[x, y] = v(xy) $. Then $v[x^{-1}, y] = v(x^{-1}y) $ is proved as for the first case. 
Since the value group is abelian, $v[x^m, y]  \ge  v( x^m y)$ holds, so we can observe
\begin{align*}
\overline{y^{-1}x^{-m} [x^m, y]} = \overline{ y^{-1}x^{-m} \sum_{\ell = 1}^m x^{m - \ell} [x, y] x^{\ell - 1}      }   =    \sum_{\ell = 1}^m \overline{y^{-1} x^{- \ell} [x, y] x^{\ell - 1}   }.
\end{align*}
For each $1 \le \ell \le m$, 
\begin{align*}
\overline{  y^{-1}x^{- \ell} [x, y] x^{\ell - 1} } &= \overline{ ( \alpha x^{-1} )^{\ell - 1}  y^{-1}x^{-1}[x, y]x^{\ell - 1}} \\
&= \overline{\alpha^{\ell - 1}x^{- \ell + 1}y^{-1}x^{-1}[x,y]  x^{\ell - 1}  } = \overline{y^{-1}x^{-1}[x,y] \alpha^{\ell - 1}}.
\end{align*}
We can change the order of $\alpha$ and $x^{-1}$ by Proposition \ref{ratdepelt} since $v(\alpha) = 0$.
The last equation follows from $v(y^{-1}x^{-1}[x,y]) = 0$ and Proposition \ref{ratdepelt}. So now we have
\begin{align*}
\overline{y^{-1}x^{-m} [x^m, y]} = 
\overline{y^{-1}x^{-1}[x,y]} \sum_{\ell = 0 }^{m - 1} \alpha^{\ell},
\end{align*}
which proves the equivalence in (2).
\endproof

\begin{trditev}\label{realrank}
Let $v$ be a valuation on a division ring $D$ such that the value group is abelian and $\overline{D} = \overline{Z(D)}$.
 Suppose the residue field is formally real and suppose $v[x, y] = v(xy) $ for some  $x, y \in D$. Then $v[x^m, y] = v(x^my)$ for all odd $m > 2$. If $v[x^m, y] > v(x^my)  $  for some even $m$, then $v[x^2, y] > v (x^2y)$ and there is no 
$a \in D$ such that $v(a^2) = v(x)$.
\end{trditev}
\proof
Suppose $v[x, y] = v(xy) $. If $m$ is odd, $\sum_{\ell = 0 }^{m - 1} \overline{\alpha}^{\ell} = 0$ does not have a solution in the residue field. By Lemma
\ref{integerpowers}, it follows that $v[x^m, y] = v(x^my)$.
Now consider the case for even $m$. If $v[x^m, y] > v (x^my)$,  then
$ \overline{\alpha} = -1$ by Lemma \ref{integerpowers} and $v[x^2, y] > v (x^2y)$. Suppose $a \in D$ satisfies
$v(a^2 ) = v(x)$. We will first show that $v[a^2, y] = v(a^2y)$. Assume that $v[a^2, y] > v(a^2y)$. Then on one hand,
\begin{align*}
\overline{a^{-2}x} = \overline{a^{-2} xy^{-1}y  }  =  \overline{y^{-1}a^{-2}xy },
\end{align*}
since  $v(a^{-2}x  ) = 0$. On the other hand,
\begin{align*}
\overline{a^{-2}x}  =  \overline{y^{-1}ya^{-2}x } = \overline{y^{-1}a^{-2}yx },
\end{align*}
where the last equation follows from $v[a^2, y] > v(a^2y)$, or, by Lemma \ref{integerpowers} equivalently, $v[a^{-2}, y] > v(a^{-2}y)$. From $\overline{x^{-1}a^{-2}(xy - yx  )} = 0$ we conclude $v[x, y] > v(xy)$, which is a contradiction.
So $v[a^2, y] = v(a^2y)$ and $v[a, y] = v(ay)$ follows from Lemma \ref{integerpowers}. 
Now we show $v[a^4, y] > v(a^4y)$. On one hand, we can write
\begin{align*}
\overline{a^{-4}x^2} = \overline{y^{-1}a^{-4}x^2y } = \overline{y^{-1}a^{-4}yx^2}
\end{align*}
since $v[x^2, y] > v(x^2y)$. On the other hand,
\begin{align*}
\overline{a^{-4}x^2} =\overline{y^{-1}ya^{-4}x^2},
\end{align*}
so we conclude $v[a^4, y] > v(a^4y)$. But by Lemma \ref{integerpowers},
$v[a^4, x] > v(a^4x)$ gives us $\overline{xax^{-1}a^{-1}} = -1$. But then, again by Lemma \ref{integerpowers}, $v[a^2, x] > v(a^2x)$. The proposition is thus proved.
\endproof

\begin{trditev} \label{realrank2}
Let $v$ be a valuation on a division ring $D$ such that  such that $\overline{D} = \overline{Z(D)}$, the value group is abelian and $2$-divisible and the residue field is 
 formally real. Suppose the value group of $v$ is of rational rank 2 and suppose there are $x, y \in D^{*}$ such that $v(x)$ and $v(y)$ are rationally independent with $v[x, y ] > v(xy)$. Then $v$ is strongly abelian.
\end{trditev}
\proof
Suppose $a, b \in D$. Since the value group is abelian, $v[a, b] \ge v(ab)$. Suppose $v[a, b] = v(ab)$. Then $v(a^{k_1}) = v(x^{-m_1}y^{- n_1})$ and
$v(b^{k_2}) = v( x^{- m_2}y^{-n_2}  ) $ for some \\ $k_i, m_i, n_i \in \Z$, $i = 1,2$. 
We conclude from Lemma 
\ref{integerpowers} that $v[a^{k_1}, b^{k_2}] = v(a^{k_1}b^{k_2}  ).$ This is immediate if $k_1$ and $k_2$ are both odd. If $k_1$ or $k_2$ is even, 
$v[a^{k_1}, b^{k_2}] = v(a^{k_1} b^{k_2}  )$ follows from the $2$-divisibility of the value group and Proposition \ref{realrank}.
Let $c := x^{m_1}y^{n_1}$ and $d := x^{m_2}y^{n_2}$. Then on the one hand,
$$\overline{a^{k_1}cb^{k_2}d} =  \overline{ca^{k_1}db^{k_2}} = \overline{dca^{k_1}b^{k_2}}$$ since $v(a^{k_1})$ and $v(c)$ are rationally dependent,  $v(b^{k_2})$ and $v(d)$ are rationally dependent and $v(b^{k_2}d) = v(a^{k_1}c)= 0$. On the other hand, $$ \overline{a^{k_1}cb^{k_2}d} = \overline{b^{k_2}da^{k_1}c} = \overline{db^{k_2}ca^{k_1}} = \overline{cdb^{k_2}a^{k_1}} = 
\overline{dcb^{k_2}a^{k_1}}.$$ Here, the last equality follows from $v[x, y] > v(xy)$ and Lemma \ref{integerpowers}. Thus we have $v(dc(a^{k_1}b^{k_2} - b^{k_2}a^{k_1}   )) > 0$, so we get  $v[a^{k_1}, b^{k_2}] > v(a^{k_1}b^{k_2}) $ which contradicts our assumption $v[a, b] = v(a, b)$. We conclude $v[a, b] > v(ab)$.
\endproof

The proof of the following proposition is the same as the proof of Proposition \ref{realrank2}.
\begin{trditev} \label{realrank22}
Let $v$ be a valuation on a division ring $D$ such that $\overline{D} = \overline{Z(D)}$, the value group is abelian and the residue field is 
formally real. Suppose the value group of $v$ is of rational rank 2 and suppose there are $x, y \in D^{*}$ such that for every $z \in D$, 
$v(z^{k}) =  v(x^{-m}y^{-n}   )$ holds for some $k, m, n \in \Z$ where $k$ is odd. Then $v$ is strongly abelian.
\end{trditev}

We will later use this result to show that all valuations $v$ on $\rweyl$ with residue field $\R$ are strongly abelian. Propositions \ref{realrank2} and \ref{realrank22} can be easily generalized to higher rational ranks of the value group. The proofs are analogous.

\begin{posledica}
Let $v$ be a valuation on a division ring $D$ such that $\overline{D} = \overline{Z(D)}$. Suppose the value group is abelian
 and $2$-divisible of rational rank $n$ and that there are $x_1, \ldots, x_n \in D$ such that $v(x_1), \ldots, v(x_n)$ are rationally independent with $v[x_i, x_j] > v(x_ix_j)$ for all $i,j$. Then $v$ is strongly abelian.
\end{posledica}

\begin{posledica}
Let $v$ be a valuation on division ring $D$ such that $\overline{D} = \overline{Z(D)}$. Suppose the value group is abelian and of rational rank $n$ and that there are $x_1, \ldots, x_n \in D$ such that for every $z \in D$, $v(z^k) = v(x_1^{m_1} \cdots x_n^{m_n})$ for some $k ,m_1, \ldots, m_n \in \Z $ with $k$ odd. 
Then $v$ is strongly abelian.
\end{posledica}

\section{Valuations on $\rweyl$}\label{rweyl}

We now describe the construction of all valuations on $\rweyl$ with residue field $\R$ that was sketched in \cite{SHT} over the ground field of $\mathbb{C}$. Since every $f \in \rweyl$ can be written as 
$\sum_{m, n \ge 0} \alpha_{m, n} x^m y^n$, the construction will be similar to the construction of all valuations on the field of rational functions $\R(x, y)$ with residue field $\R$ (examples of constructions of such valuations can be found in \cite{KUHL} or \cite{MAC}), but with some additional constraints arising from the fact that the generators $x, y \in \rweyl$ satisfy $[y, x] = 1$. We first note that it follows from Theorem 5.3 of \cite{MA2} that the value group of any valuation on $\rweyl$ is commutative. Also, since every valuation $v$ on $\rweyl$ can be uniquely extended to its quotient division ring $\mathcal{D}_1(\R)$, our construction will take place in the quotient ring as we will use inverses.

To construct a valuation $v$ trivial on $\R$ with residue field $\R$, we compare $v(x)$ and $v(y)$. It is easy to show, as it was done in \cite{MA2}, that $v(xy) = v(yx) < 0$, so $v(x)$ or $v(y)$ will be less than zero. Without loss of generality, we can set $v(x) = -1 \in \Q$ and compare it to $v(y)$. If $v(y) \not \in \Q$, then we get $$v(\sum_{m, n \ge 0} \alpha_{m, n} x^m y^n) = \min_{m, n} \{ m v(x) + n v(y) \}$$ for all elements of $\rweyl$. Otherwise, $v(y) = \frac{m_1}{n_1} \in \Q$. It follows that $\overline{x^{m_1} y^{n_1} } = \beta_1 \in \R$, so
$v(x^{m_1} y^{n_1} - \beta_1) > 0$. Set $$\omega_1 := x^{m_1} y^{n_1} - \beta_1$$ and as before compare $v(\omega_1)$ to $v(x)$ in terms of rational dependence. 
If $v(\omega_1) = \frac{m_2}{n_2} \in \Q$, then $\overline{x^{m_2}\omega_1^{n_2}  } = \beta_2$ for some $\beta_2 \in \R$. Hence $$\omega_2 := x^{m_2}\omega_1^{n_2}  - \beta_2$$ also has value greater than zero.
We continue this procedure. If we additionally define  $\omega_{-1} = x$ and $\omega_{0} = y$, we thus get a sequence $$(\omega_i )_{i \ge -1}, \omega_i \in \rweyl$$ which ends with $\omega_n$ for some $n \in \N$ if $v(\omega_n) \not\in \Q$ or is infinite otherwise. 

By the end of this section, we will prove a necessary and sufficient condition for the possibility to extend $v$  from $(\omega_i)_{i \ge -1}$ to a valuation on $\rweyl$ with residue field $\R$. Every such extension from $ (\omega_i)_{i \ge -1} $ to $\rweyl$ will be uniquely determined.
We will also show that every valuation on $\rweyl$ with residue field $\R$ is strongly abelian. 

\subsection{Properties of the sequence $(\omega_i)_{i \ge -1}$ associated to a valuation on $\rweyl$}
Thorough this subsection let $v$ be a valuation on $\rweyl$.
\begin{lema} \label{xkom}
Suppose $( \omega_i  )_{i \ge -1} \subseteq \rweyl$ is a sequence as described above, with $\omega_{-1} =x$, $\omega_0 = y$,
$\omega_i = x^{m_i} \omega_{i - 1}^{n_i} - \beta_i$ for all $i \ge 0$. 
Then $[\omega_i, x] $ equals
\begin{align*}
x^{m_i} \sum_{\ell_i = 1}^{n_i} \omega_{i - 1}^{n_i - \ell_i} \bigg(  x^{m_{i - 1}}   \sum_{\ell_{i - 1} = 1}^{n_{i - 1}} \omega_{i-2}^{n_{i- 1} - \ell_{i - 1}} \bigg( \cdots \bigg(x^{m_2} \sum_{\ell_2 = 1}^{n_2 } \omega_1^{n_2 - \ell_2} n_1 x^{m_1} y^{n_1 - 1} \omega_1^{\ell_2} \bigg)\omega_2^{\ell_3} \bigg) \cdots \bigg)\omega_{i-1}^{\ell_i}
\end{align*}
for each $i \ge 1$. 
\end{lema}
\proof
We prove the lemma by induction on $i$. If $i = 1$,
\begin{align*}
[\omega_1, x] = [x^{m_1}y^{n_1} - \beta_1, x] = x^{m_1} \sum_{\ell_1 = 1}^{n_1} y^{n_1 - \ell_1} [y, x]y^{\ell_1 - 1} = n_1x^{m_1}y^{n_1 - 1}.
\end{align*}
Now suppose that the equality holds for $[\omega_{i}, x]$. Then we have
\begin{align*}
[ \omega_{i + 1}, x   ] &= [x^{m_{i + 1}}\omega_{i}^{n_{i + 1}} - \beta_{i + 1}, x] = x^{m_{i + 1}}[ \omega_i^{n_{i + 1}}, x ] \\
 &=  x^{m_{i + 1}}  \sum_{\ell_{i + 1} = 1}^{n_{i + 1}} \omega_i^{n_{i + 1} - \ell_{i  + 1}  } [ \omega_i, x    ] \omega_i^{\ell_{i + 1} - 1} 
\end{align*}
We can then proceed by the induction hypothesis.
\endproof

Before proving the next lemma, we define an equivalence relation between nonzero elements of $\rweyl$ that have the same $v$-value, but their difference does not. For any $a, b \in \rweyl \setminus \{0 \}$, we write $a \sim b$ if $v(a) = v(b) < v(a - b)$. This is also a congurence relation, as $ac \sim bc$ and $ca \sim cb$ holds for all $a, b,c \in \rweyl \setminus \{0 \}$ with $a \sim b$.

\begin{lema} \label{kompogoj} 
Suppose $v$ is a valuation on $\rweyl$ with residue field $\R$ and suppose $( \omega_i  )_{i \ge -1}$ is a sequence such that $\omega_{-1 } = x$, $\omega_0 = y$,
$\omega_i = x^{m_i} \omega_{i - 1}^{n_i} - \beta_i$,  $v(x) = -1$ and  $v(\omega_i) = \frac{m_{i + 1}}{ n_{i + 1} }$ for all $i \ge 0$ up to either some $n \ge 0$ in which case $v( \omega_n  ) \not\in \Q$,     
or up to infinity. Then $v[\omega_j, \omega_i ] > v( \omega_i \omega_j )$ for all $i, j \le k$  if and only if $v( \Pi_{ \ell = -1}^k  \omega_{\ell}   ) < 0$, where $k \le n$ in case $v(\omega_n ) \not\in \Q$ for some $n \ge 0$.

If any and hence both sides of the equivalence hold, then $$v[ \omega_j, \omega_i  ] = - v(xy \omega_1 \cdots \omega_{i-1} \omega_{i + 1} \cdots \omega_{j - 1})$$ for all $i < j \le k$.
\end{lema}

\proof
Suppose $v$ is a valuation on $\rweyl$ and $( \omega_i )_{i \ge - 1  }$ is a sequence as described in the lemma. So  $v(\omega_i) \in \Q$ either for all $i \ge 0$ or for all $0 \le i < n$ for some $n \ge 0$ and $v(\omega_n) \not\in \Q$. It follows from Proposition \ref{ratdepelt} that 
$v[\omega_i, \omega_j] > v( \omega_i \omega_j)$ for all $i,j  < n$ since $v(\omega_i)$ and $v(\omega_j)$ are rationally dependent. We shall use this fact to evaluate $v[\omega_i, \omega_j]$ for all $i, j \le k \le n$.
 It follows from Lemma \ref{xkom} that $[\omega_k, x]$ is a sum of 
products $P$, all equal to $y^{n_1 - 1}x^{m_1} \omega_1^{n_2 - 1} x^{m_2} \cdots \omega_{k - 1}^{n_k - 1}x^{m_k} $ up to the order of factors. Since 
$v[ \omega_i , \omega_j] > v( \omega_i \omega_j )$ for all $i,j \le k-1$, 

$$P \sim y^{n_1 - 1}x^{m_1} \omega_1^{n_2 - 1} x^{m_2} \cdots \omega_{k - 1}^{n_k - 1}x^{m_k} $$
holds for every product $P$ of the sum. Since
\begin{align*} 
v( y^{n_1 - 1}x^{m_1} \omega_1^{n_2 - 1} x^{m_2} \cdots \omega_{k - 1}^{n_k - 1}x^{m_k}    ) + v(y \omega_1 \cdots \omega_{k - 1}) = \sum_{i = 1}^{k}v(  x^{m_i}\omega_{i-1}^{n_i}) = 0, 
\end{align*}
we can conclude $v[ \omega_k, x ] = v( y^{n_1 - 1}x^{m_1} \omega_1^{n_2 - 1} x^{m_2} \cdots \omega_{k - 1}^{n_k - 1}x^{m_k}    ) =  - v(y \omega_1 \cdots \omega_{k - 1})$. It follows that $v[x, \omega_k] > v(x\omega_k)$ if and only if $v(xy \omega_1 \cdots \omega_k  ) < 0$.

We will now prove that $v[ \omega_{i + k}, \omega_i  ] = - v(xy \omega_1 \cdots \omega_{i-1} \omega_{i + 1} \cdots \omega_{i + k - 1})$ by induction on $k$,  $1 \le k \le n - i$. 
It will then follow that $v[\omega_i, \omega_j ] > v( \omega_i \omega_j )$ for all $i, j \le n$ if and only if $v( \Pi_{\ell = -1}^n  \omega_{\ell}   ) < 0$.
If $k  = 1$, then 
\begin{align*}
[\omega_{i + 1}, \omega_i ] &= [ x^{m_{i + 1}} \omega_i^{n_{i + 1}} - \beta_{i + 1} , \omega_i ] = [x^{m_{i + 1}}, \omega_i   ] \omega_i^{n_{i + 1}} \\
 &= ( \sum_{\ell = 1}^{m_{i  + 1}} x^{m_{i + 1} - \ell  } [ x, \omega_i ] x^{\ell -1}  ) \omega_i^{n_{i + 1}},
\end{align*}
and since $[x, \omega_i]$ is a sum of products all equal to $y^{n_1 - 1}x^{m_1} \omega_1^{n_2 - 1} x^{m_2} \cdots \omega_{i - 1}^{n_i - 1}x^{m_i} $    up to the order of factors, we can, using $v[ \omega_i, \omega_j] > v( \omega_i \omega_j  )$ for $j < i$, deduce that $v[\omega_{i + 1}, \omega_i] = 
- v(xy\omega_1 \cdots \omega_{i - 1}  )$ just like we did when evaluating $v[\omega_i, x]$. For $k > 1$, we have
\begin{align*}
[ \omega_{i + k}, \omega_i  ] &= [x^{m_{i + k}} \omega_{i + k - 1}^{n_{i + k}} - \beta_{i + k}, \omega_i  ] = [x^{m_{i + k}}, \omega_i  ] \omega_{i + k - 1}^{n_{i + k}} + x^{m_{i + k}}[  \omega_{i + k - 1}^{n_{i + k}}   ,  \omega_i   ] \\
&= (  \sum_{ \ell = 1}^{m_{i + k}} x^{m_{i + k }  -  \ell} [x, \omega_i] x^{\ell - 1}       )\omega_{i + k - 1}^{n_{i + k}} + x^{m_{i + k}}( \sum_{ \ell  = 1}^{n_{i + k}} 
\omega_{i+k - 1}^{n_{i + k } - \ell} [ \omega_{i + k - 1}, \omega_i  ] \omega_{i + k - 1}^{\ell}   ) \\
&\sim m_{i + k} x^{m_{i + k} - 1}\omega_{i + k - 1}^{n_{i + k}} [x, \omega_i] + n_{i + k} x^{m_{i + k}} 
\omega_{i + k - 1}^{n_{i + k} - 1} [ \omega_{i + k - 1}, \omega_i  ] 
\end{align*}
and using both Lemma \ref{xkom} and induction on $k$, we see that the first sum has $v$-value equal to $ - v( xy \omega_1 \cdots \omega_{i - 1}   ) $ and the second has $v$-value equal to
$- v(xy\omega_1 \cdots \omega_{i + k - 1} ) $. Since the latter is smaller, it is equal to $v[\omega_{i+k}, \omega_i]$. This proves the lemma.
\endproof

It follows that if $v$ can be extended from $( \omega_i)_{i  \ge -1 }$ to a valuation on $\rweyl$, $ \sum_{i \ge -1}^k v( \omega_i  )$ must be strictly less than $0$ for all $k \le n$ in case $v(\omega_n) \not\in \Q$ for some $n$, and for all $k \ge 0$ if $v(\omega_i) \in \Q$ for all $i \in \N$. We will now describe a necessary condition for the residue field to be $\R$ and then proceed to show that if both conditions are fulfilled, $v$ can be extended from $( \omega_i   )_{i \ge -1}$  to a valuation on $\rweyl$ with residue field $\R$. 

To ensure that the residue field is $\R$, it is obviously necessary that $\overline{ \omega_{i - 1}^{k_i} \omega_{j - 1}^{k_j}   }  \in \R $ holds for all $k_i, k_j \in \Z$ with $k_i \frac{m_i}{n_i} + k_j \frac{m_j}{n_j} = 0$. For given $i , j \ge 0$, all  solutions $(k_i, k_j) \in \Z^2$ to the diophantine equation
\begin{align}\label{diofantska}
k_i m_jn_i + k_j m_in_j = 0
\end{align}
are integer multiples of the pair $(K_{i,j},-  K_{j, i})$ with 
\begin{align*}
K_{i, j} &= \frac{m_jn_i}{d_{i,j}}, \\
K_{j, i} &=  \frac{m_in_j}{d_{i,j}}, \\
\text{where }d_{i,j} &= \gcd \{ m_jn_i, m_in_j \}.
\end{align*}
So for all $k_i, k_j \in \Z$ with $k_i \frac{m_i}{n_i} + k_j \frac{m_j}{n_j} = 0$, we can write
\begin{align*}
\overline{ \omega_{i - 1}^{k_i} \omega_{j - 1}^{k_j}   } = \overline{ \omega_{i - 1}^{nK_{i,j}} \omega_{j - 1}^{-nK_{j,i}}   } = 
( \overline{\omega_{i - 1}^{K_{i,j}} \omega_{j - 1}^{ -  K_{j, i}}}   )^n 
\end{align*}
for some $n \in \Z$, where we used Proposition \ref{ratdepelt} in the second equality. So for every $k_i, k_j \in \Z$ satisfying \ref{diofantska},    $\overline{ \omega_{i - 1}^{k_i} \omega_{j - 1}^{k_j}   }$ is uniquely determined by $ \overline{\omega_{i - 1}^{K_{i,j}} \omega_{j - 1}^{ - K_{j, i}}} $. 
For each $i, j \ge 0$, we define $\alpha_{i, j} = \overline{\omega_{i - 1}^{K_{i,j}} \omega_{j - 1}^{ - K_{j, i}}}$. We immediately see that $\alpha_{j ,i} = \alpha_{i, j }^{ - 1}$ and hence $\alpha_{i,i} = 1$ for all $i, j \ge 0$. As
\begin{align*}
\alpha_{i,j }^{d_{i,j}} =  
\overline{ \omega_{i - 1}^{ K_{i,j} d_{i,j} }  \omega_{j - 1}^{- K_{j,i} d_{i,j} }   }   =  \overline{\omega_{ i - 1 }^{ m_jn_i  } \omega_{j - 1}^{ - m_in_j }  } = \beta_i^{m_j} \beta_j^{- m_i}
\end{align*}
for all $i, j \ge 0$, $\alpha_{i,j}$ is one of the possible $d_{i,j}$-th roots for $\beta_i^{m_j} \beta_j^{- m_i}.$
If $v$ is a valuation on $\mathcal{D}_1( \R )$ with residue field $\R$, $\alpha_{i, j}$ must be real for all $i, j \ge 0$. For every $i,j \ge 0$ with even $d_{i,j}$, this means that $\beta_i^{m_j} \beta_j^{- m_i} > 0$ must hold.
In the next lemma, we present a necessary condition on the sequence $( \beta_i )_{i \ge 1}$ so that $\alpha_{i, j} \in \R$ can be chosen for all $i,j$. We also prove that if $n_i$ is odd, $\alpha_{i,j}$ is uniquely determined for all $j \ge 0$.

\begin{lema}\label{dvojice}
Let $v$ be a valuation on $\mathcal{D}_1(\R)$ as in Lemma  \ref{kompogoj}. Then the following holds:
\begin{enumerate}
\item If $n_i$ is odd, there is a unique possible choice for $\alpha_{i,j} \in \R$ for all $j \ge 0$.
\item Only if $\sgn(\beta_i)$ is constant on the set of all $i \ge 0$ for which $n_i$ is even can we choose $\alpha_{i ,j } \in \R$ for all $i, j \ge 0$.
\end{enumerate}
\end{lema}
\proof
Suppose $n_i$ is odd. Then for any $j \ge 0$, let $\tilde{d}_{i,j}$ be the highest odd number dividing $d_{i,j}$. Since $n_i$ is odd,
$ \ell_1 := \frac{\tilde{d}_{i,j} m_j}{ d_{i,j}  } \in \Z$. If $n_j$ is odd as well,  $\ell_2 := \frac{\tilde{d}_{i,j} m_i}{ d_{i,j}} \in \Z$ holds too.
If $n_j$ is even, $m_j$ is odd, so $d_{i, j} = \tilde{d_{i, j}}$ as $d_{i, j}$ 
divides $ m_jn_i $. In both cases, $\ell_2  \in \Z$ holds. Then for $\ell := \ell_1 m_i = \ell_2m_j =  \frac{\tilde{d}_{i,j} m_im_j}{ d_{i,j}  }$ we can evaluate
\begin{align*}
\alpha_{i,j}^{ \tilde{d}_{i,j} } = \overline{\omega_{i - 1}^{ K_{i,j}\tilde{d}_{i,j}  } \omega_{j - 1}^{ - K_{j,i}\tilde{d}_{i,j}} } = 
\overline{x^{ \ell  }   x^{-  \ell  } \omega_{i - 1}^{ \ell_1 n_i  } \omega_{j - 1}^{ - \ell_2 n_j}} = 
\overline{ (x^{m_i}\omega_{i-1}^{n_i} )^{\ell_1} (x^{m_j}\omega_{j-1}^{n_j })^{- \ell_2}   } = \beta_i^{\ell_1}\beta_j^{- \ell_2},
\end{align*}
and since $\tilde{d}_{i,j}$ is odd, $\alpha_{i,j} \in \R$ is uniquely determined. The first point of the lemma is thus proven.

To prove the second point of the lemma, suppose $i, j \ge 0$ are such that $n_i$ and $n_j$ are both even. As a consequence, both $m_i$ and $m_j$ are odd while $d_{i,j}$ is even. So, provided $\alpha_{i, j } \in \R$, we compute
\begin{align*}
1 = \sgn( \alpha_{i, j}^{d_{i,j}} ) = \sgn( \beta_i^{m_j} \beta_j^{-m_i} ) = \sgn( \beta_i \beta_j ),
\end{align*}
which proves the second part of the lemma.
\endproof
For even $d_{i,j}$ we have seemingly two choices for $\alpha_{i,j} \in \R$ -- a positive and a negative one. We will show that in most cases, we cannot choose $\sgn(\alpha_{i,j})$ for all $i, j \ge 0$ independently of each other. 

Before that, we observe that for any $i, j \ge 0$, at most one of $K_{i,j}$ and $K_{j, i}$ is even. In fact, if at most one of $n_i$ and $n_j$ is odd, $K_{i,j}$ is odd if and only if $n_i$ is divisible by the greatest power of two that divides $n_j$. For each $i \ge 1$, let $2^{h_i}$ be the biggest power of two that divides $n_i$. Define also $m_0 = 1$ and  $n_0 = -1$.

\begin{trditev} \label{dvojicedet}
Let $v$ be a valuation on $\mathcal{D}_1(\R)$ associated to a sequence $( \omega_i)_{i \ge -1}$ with $v(\omega_{-1} ) = v(x) = -1$, $v(\omega_{i - 1}) = \frac{m_{i}}{n_{i}}$ with $\gcd(m_{i },  n_{i}) = 1$  and $\overline{ x^{m_i} \omega_{i - 1}^{n_i} } = \beta_i \in \R$ for each $i \ge 1$. Suppose $\sgn(\beta_i)$ is constant on the set of all $i \ge 0$ for which $n_i$ is even. Suppose $ \alpha_{i,j} \in \R$ is determined for all $i, j \ge 0$. Then 
$\overline{\Pi_{i = 0}^r \omega_{i-1}^{k_i}} \in \R$ is uniquely determined for each 
set of integers $k_0, k_1, \ldots, k_r \in \Z$ with $\sum_{i = 0}^r k_i \frac{m_i}{n_i} = 0$ if and only if for each $a, b ,c \ge 0$,
$\alpha_{a,b}\alpha_{a,c}\alpha_{b,c} > 0$ whenever $h_a = h_b \le h_c$ holds. 
\end{trditev}
\proof
To prove the necessity of the condition, suppose $a, b, c \ge 0$ are such that $h_a = h_b \le h_c$ holds. Suppose $K_{a,b}$ and $K_{b,c}$ are both odd. Choose
 $k_a, k_b, k_c \in \Z \setminus \{ 0 \}$ such that $k_a \frac{m_a}{n_a  } +k_b \frac{m_b}{n_b  } + k_c \frac{m_c}{n_c  } = 0$ and that $k_a$ and $k_b$ are odd while $k_c$ is even. Then on one hand,
\begin{align*}
 \sgn ( \overline{ \omega_{a - 1}^{k_a} \omega_{b - 1}^{k_b} \omega_{c - 1}^{k_c} }    )  &=  
\sgn ( \overline{ \omega_{a - 1}^{k_a} \omega_{b - 1}^{k_b} \omega_{c - 1}^{k_c} } ^{K_{a,b}}   ) = \sgn(   \overline{ \omega_{a - 1}^{k_a} \omega_{b - 1}^{k_b} \omega_{c - 1}^{k_c} } ^{K_{a,b}} \overline{ \omega_{b-1}^{- k_a K_{b, a}}  \omega_{b-1}^{ k_a K_{b, a}} }   )  \\
&= \sgn(   \alpha_{a,b}^{k_a} \overline{ \omega_{ b - 1 }^{  {\ell}_b  } \omega_{ c - 1 }^{  {\ell}_c }   } )
\end{align*}
with 
\begin{align*}
\ell_b &= k_bK_{a,b} + k_a K_{b,a},  \\ 
\ell_c &= k_c K_{a,b}.
\end{align*}
As $v(\omega_{ b - 1 }^{  {\ell}_b  } \omega_{ c - 1 }^{  {\ell}_c }) = 0$, $(\ell_b, \ell_c) = \ell(K_{b,c}, - K_{c, b})$ for some $\ell \in \Z$ with $- \ell K_{c,b} = \ell_c = k_c K_{a,b} $. So we can conclude 
\begin{align*}
\sgn ( \overline{ \omega_{a - 1}^{k_a} \omega_{b - 1}^{k_b} \omega_{c - 1}^{k_c} }    )  = \sgn(\alpha_{a,b}^{k_a}\alpha_{b,c}^{\ell}).
\end{align*}
On the other hand, we see by analogous computations  that 
\begin{align*}
 \sgn ( \overline{ \omega_{a - 1}^{k_a} \omega_{b - 1}^{k_b} \omega_{c - 1}^{k_c} }    )  &=  \sgn ( \overline{ \omega_{a - 1}^{k_a} \omega_{b - 1}^{k_b} \omega_{c - 1}^{k_c} } ^{K_{a,c}}   ) = \sgn(   \overline{ \omega_{a - 1}^{k_a} \omega_{b - 1}^{k_b} \omega_{c - 1}^{k_c} } ^{K_{a,c}} \overline{ \omega_{c-1}^{- k_a K_{c, a}}  \omega_{c-1}^{ k_a K_{c, a}} }   )  \\
&= \sgn(   \alpha_{a,c}^{k_a} \overline{ \omega_{ b - 1 }^{  {\ell'}_b  } \omega_{ c - 1 }^{  {\ell'}_c }   } ) = 
\sgn(   \alpha_{a,c}^{k_a}   \alpha_{b,c}^{\ell'}  )
\end{align*}
for some $\ell'_b, \ell'_c, \ell' \in \Z$ with 
\begin{align*}
\ell'_b &= k_b K_{a,c} = \ell' K_{b,c}  ,  \\ 
\ell'_c &= k_c K_{a,c} + k_a K_{c,a}. 
\end{align*}
We have chosen $k_a, k_b$ odd and $k_c$ even. In this case, the greatest power of two that divides $k_c$ is $2^{h_c - h_a + 1}$. On the other hand, the greatest power of two that divides $K_{c,a}$ and $K_{c,b}$ is $2^{h_c - h_a}$. We can thus conclude from $\ell K_{c,b} = - k_c K_{a,b} $ and $\ell' K_{b,c} =  k_b K_{a,c}$ that $\ell$ is even while $\ell'$ is odd since $K_{a,b}, K_{a,c}$ and $K_{b,c}$ are all odd. So we see that 
\begin{align*}
 \sgn ( \overline{ \omega_{a - 1}^{k_a} \omega_{b - 1}^{k_b} \omega_{c - 1}^{k_c} }    )  = \sgn(\alpha_{a,b})
= \sgn( \alpha_{a,c} \alpha_{b,c}   ),
\end{align*}
which proves the necessity of the condition.

Now suppose $\sgn ( \alpha_{a,b} \alpha_{a,c}\alpha_{b,c}    ) = 1$ for all $a, b, c \ge 0$ with 
$h_a = h_b \le h_c$. Let $K := (k_{0}, \ldots, k_r   ) \in \Z^{r + 1} $ be such that 
$\sum_{i = 0}^{r}k_i \frac{m_i}{n_i} = 0$. Let $\supp K := \{ i \mid k_i \neq 0 \}$ and $n := |\supp K |$. 
We prove that $ \overline{ \Pi_{i = 0}^r \omega_{i - 1}^{k_i}    } \in \R $ is uniquely determined by induction on $n \ge 2$. We first suppose $0\not \in \supp K.$ We will deal with the case 
$0 \in \supp K$ at the end of our proof.

If $n = 2$, then $\Pi_{i =0}^r \omega_{i - 1}^{k_i} = \omega_{i - 1}^{k_i}\omega_{j - 1}^{k_j} $ for some $i,j > 0$, and its value in the residue field is a power of $\alpha_{i,j}$.

Now suppose $n > 2$. Take two distinct $a, b \in \supp K$. As at least one of $K_{a,b}$ and $K_{b, a}$ is odd, so suppose $K_{a, b}$ is odd. Then
\begin{align*}
\overline{ \Pi_{i = 1}^r \omega_{i - 1}^{k_i}    }^{K_{a,b}} = \overline{ \Pi_{i = 1}^r \omega_{i - 1}^{k_i}    }^{K_{a,b}} \overline{\omega_{b - 1}^{ -k_a K_{b,a}  } \omega_{b - 1}^{ k_a K_{b,a}    }} = 
\overline{  \Pi_{i = 1}^r \omega_{i - 1}^{\ell_{i, 1}}    }   \alpha_{a,b}^{k_a}
\end{align*}
with $\ell_{a, 1} = 0, \ell_{b, 1} = k_a K_{a,b} + k_b K_{b, a} $ and $\ell_{i, 1} = k_i K_{a, b}$ for $i \neq a, b$.
Since $| \{ i \mid \ell_{i, 1} \neq 0  \} |$ is strictly smaller than $n$, 
$ \overline{  \Pi_{i = 0}^r \omega_{i - 1}^{\ell_{i, 1}}    } \in \R$ is uniquely determined by the induction hypothesis. 
So we have determined $\overline{ \Pi_{i = 0}^r \omega_{i - 1}^{k_i}    }^{K_{a,b}} \in \R $. As $K_{a,b}$ is odd, 
$\overline{ \Pi_{i = 0}^r \omega_{i - 1}^{k_i}    } \in \R$ is determined as well.

We now need to show that in this way, $\overline{ \Pi_{i = 1}^r \omega_{i - 1}^{k_i}  }$ is uniquely determined, that is, if we choose another $a', b' \in \supp K$ instead of $a, b$, we get the same value for 
$\overline{ \Pi_{i = 1}^r \omega_{i - 1}^{k_i}    } \in \R$. We will show this by choosing $c \in \supp K \setminus \{ a, b \}$ and proving that the evaluated value of $\overline{ \Pi_{i = 1}^r \omega_{i - 1}^{k_i}    }$ is the same whether we factor a power of $\alpha_{a,b}$ as above, or $\alpha_{a, c}$ or $\alpha_{b,c}$ instead. By transitivity of the equality relation, this will imply that the obtained value of 
$\overline{ \Pi_{i = 1}^r \omega_{i - 1}^{k_i}  }$ is independent of the choice of $a, b \in \supp K$. Suppose without loss of generality that $h_a \le h_b \le h_c$ and that $K_{a,b}, K_{a,c}$ and $K_{b,c}$ are odd.
Above, we have evaluated
\begin{align*}
\overline{ \Pi_{i = 1}^r \omega_{i - 1}^{k_i}    }^{K_{a,b}} = \overline{  \Pi_{i = 1}^r \omega_{i - 1}^{\ell_{i, 1}}    }   \alpha_{a,b}^{k_a}
\end{align*}
with $\ell_{a, 1} = 0, \ell_{b, 1} = k_a K_{a,b} + k_b K_{b, a} $ and $\ell_{i, 1} = k_i K_{a, b}$ for $i \neq a, b$. We proceed by evaluating, in the same way as before,
\begin{align*}
\overline{  \Pi_{i = 1}^r \omega_{i - 1}^{\ell_{i, 1}}    }^{K_{b, c} } = \overline{  \Pi_{i = 1}^r \omega_{i - 1}^{p_{i, 1}}    } \alpha_{b,c}^{ \ell_{b, 1}  }
\end{align*}
with $p_{a,1} = p_{b, 1} = 0$, $p_{c, 1} = \ell_{c, 1}K_{b,c} + \ell_{b,1} K_{c, b} $ and $p_{i, 1} = \ell_{i, 1}K_{b,c}$ for $i \neq a, b, c$.
So
\begin{align} \label{k1}
\overline{ \Pi_{i = 1}^r \omega_{i - 1}^{k_i}    }^{K_{a,b}K_{b,c}} = \overline{\Pi_{i = 1}^r \omega_{i - 1}^{p_{i, 1}} } \alpha_{a,b}^{k_a K_{b,c}} \alpha_{b,c}^{ \ell_{b, 1}  } 
\end{align}
with $\ell_{i,1}$ and $p_{i,1}$ for all $0 \le i \le r$ as above. In particular, we see that for $i \neq a, b, c$, $p_{i, 1} = \ell_{i, 1}K_{b,c} = k_i K_{a,b} K_{b,c}$. 
Similarly, we can compute 
\begin{align} \label{k2}
\overline{ \Pi_{i = 1}^r \omega_{i - 1}^{k_i}    }^{K_{a,c}K_{b,c}} = (\overline{  \Pi_{i = 1}^r \omega_{i - 1}^{\ell_{i, 2}}    } \alpha_{a,c}^{k_a})^{K_{b,c}} = 
\overline{  \Pi_{i = 1}^r \omega_{i - 1}^{p_{i, 2}}  }  \alpha_{a,c}^{k_aK_{b,c}} \alpha_{b,c}^{\ell_{b,2}}
\end{align}
with 
\begin{enumerate}
\item $\ell_{a, 2} = 0, \ell_{c, 2} = k_a K_{a,c} + k_c K_{c, a} $ $\ell_{i, 1} = k_i K_{a, c}$ for $i \neq a, c$, and      \item $p_{a,2} = p_{b, 2} = 0, p_{c, 2} = \ell_{b,2}K_{c,b} + \ell_{c,2} K_{b,c}$, $p_{i, 2} = k_i K_{a,c}K_{b,c}$ for $i \neq a, b, c$. 
\end{enumerate}
Let $N := K_{a,b}K_{a,c}K_{b,c}$. On one hand, we see from \ref{k1} that
\begin{align} \label{produkt1}
\overline{ \Pi_{i = 1}^r \omega_{i - 1}^{k_i}    }^N = \overline{\Pi_{i = 1}^r \omega_{i - 1}^{p_{i, 1}} }^{K_{a,c}} \alpha_{a,b}^{k_a  K_{b,c}K_{a,c}} \alpha_{b,c}^{ \ell_{b, 1}K_{a,c}  },
\end{align}
and on the other hand, we see from \ref{k2} that
\begin{align} \label{produkt2}
\overline{ \Pi_{i = 1}^r \omega_{i - 1}^{k_i}    }^N = \overline{\Pi_{i = 1}^r \omega_{i - 1}^{p_{i, 2}} }^{K_{a,b}} \alpha_{a,c}^{k_a K_{b,c} K_{a,b}} \alpha_{b,c}^{ \ell_{b, 2}K_{a,b}  }.
\end{align}
We need to show that in both equations, we get the same value. We first see that for all $i \neq c$, $p_{i, 1} K_{a,c} = p_{i, 2} K_{a,b}$. So, given that
$$\sum_{i = 1}^r p_{i, 1} \frac{m_i}{n_i} = \sum_{i = 1}^r p_{i, 2} \frac{m_i}{n_i} = 0, $$
we can see  $p_{c, 1} K_{a,c} = p_{c, 2} K_{a,b}$ holds as well, and thus we conclude
$$\overline{\Pi_{i = 1}^r \omega_{i - 1}^{p_{i, 1}  }}^{K_{a,c}} = \overline{\Pi_{i = 1}^r \omega_{i - 1}^{p_{i, 2}  } }^{K_{a,b}}.$$ 
It then follows that
$$\alpha_{a,b}^{k_a K_{a,c} K_{b,c}} \alpha_{b,c}^{ \ell_{b, 1}K_{a,c}  } = \alpha_{a,c}^{k_a K_{b,c} K_{a,b}} \alpha_{b,c}^{ \ell_{b, 2}K_{a,b}  },$$
since both sides of the equation are equal to $\overline{\omega_{a - 1}^{N k_a} \omega_{b - 1}^{N k_b} \omega_{c - 1}^{N k_c - p_{c, 1} K_{a,c}}}$ and the signs of $\alpha_{a,b}, \alpha_{a,c}$ and $\alpha_{b, c}$ were chosen so that the signs of both sides of the equality match. We conclude that the value of $\overline{ \Pi_{i = 1}^r \omega_{i - 1}^{k_i}    }^N$ is the same in both \ref{produkt1} and 
\ref{produkt2}. As $N$ is odd (since $K_{a,b}$, $K_{a,c}$ and $K_{b,c}$ are all odd), we conclude that $\overline{ \Pi_{i = 1}^r \omega_{i - 1}^{k_i}    }$ is the same whether we factor a power of $\alpha_{a,b}$ or $\alpha_{a, c}$. If we factored a power of $\alpha_{b,c}$,
 we would, as similar computations as above would show, get the same value for $\overline{ \Pi_{i = 1}^r \omega_{i - 1}^{k_i}    }$.

We have now shown that if the condition of the proposition is fulfilled, $\overline{ \Pi_{i = 1}^r \omega_{i - 1}^{k_i}    }$ is uniquely determined for all $k_1, \ldots, k_r \in \Z$ with $\sum_{i = 1}^r k_i \frac{m_i}{n_i} = 0$. 

Now we consider the case $0 \in \supp K$. Let $K = (k_0, \ldots, k_r) \in \Z^{r + 1}$ be such that $k_0 \neq 0$ and 
$\sum_{i = 0}^{r} k_i \frac{m_i}{n_i} = 0$. Let $N := \gcd \{ n_i \mid i \in \supp K  \}$. If $N$ is odd, i.e., if $n_i $ is odd for every $i \in \supp K$, then $\overline{\Pi_{i = 0}^{r} \omega_{i-1}^{k_i}   }$ is uniquely determined. This is because 
\begin{align*}
{\overline{\Pi_{i = 0}^{r} \omega_{i-1}^{k_i}   }^N} = 
(\overline{\Pi_{i= 0}^r x^ {  k_{i } \frac{m_i}{n_i}}  {\Pi_{i = 0}^{r} \omega_{i-1}^{k_i}   }})^N    = 
\overline{\Pi_{i = 1} (  x^{m_i} \omega_{i - 1}^{ni}  )^{k_ic_i}} = \Pi_{i = 1}^{r} \beta_i^{k_ic_i}
\end{align*}
where $c_i := \frac{N}{n_i}$ for each $i \le i \le r$. We thus conclude
$ \overline{\Pi_{i = 0}^{r} \omega_{i-1}^{k_i}   } \in \R $ is the uniquely determined $N$-the real root of $\Pi_{i = 1}^{r} \beta_i^{k_ic_i}$.
Now suppose $n_j$ is even for some $j \in \supp K$. Then $m_j$ must be odd since $\gcd (m_j, n_j ) = 1$.
Let $k'_j := k_j - k_0n_j$ and $k'_i := k_i$ for all $i \in \supp K \setminus \{0, j \}$. Then $\sum_{i = 1}^{r} k'_i \frac{m_i}{n_i} = 0$ and
\begin{align*}
{\overline{ \Pi_{i = 0}^r \omega_{i - 1}^{k_i}    }}^{m_j  } = (\overline{x^{m_j} \omega^{n_j}   })^{k_0}
(\overline{\Pi_{i = 1}^{r} \omega_{i-1}^{k'_i}   })^{m_j} = \beta_j^{k_0}(\overline{\Pi_{i = 1}^{r} \omega_{i-1}^{k'_i}   })^{m_j}.
\end{align*}
We evaluate $(\overline{\Pi_{i = 1}^{r} \omega_{i-1}^{k'_i}   })^{m_j}$ as above $k_0  = 0$ and conclude that 
${\overline{ \Pi_{i = 0}^r \omega_{i - 1}^{k_i}    }} \in \R $ is the unique $m_j$-th real root of 
$ \beta_j^{k_0}(\overline{\Pi_{i = 1}^{r} \omega_{i-1}^{k'_i}   })^{m_j}$.

This concludes the proof of our proposition.
\endproof

In Lemma \ref{babyconstruct}, we suppose that $v$ is a valuation on $\mathcal{D}_1(\R)$ extended from $ ( \omega_i)_{i \ge - 1}$ to 
$\mathcal{D}_1( \R)$ and compute the value of certain elements of $\mathcal{D}_1(\R)$ in this case.

\begin{lema} \label{potence}
Let $D$ be a division ring endowed with a valuation $v$ with an abelian value group and a commutative residue field with characteristic zero. Let $a, b \in D$ be such that $a \sim b$, $v(a) = v(b) = 0$ and $v(ab) < v[a, b]$. Then $v(a^n - b^n) = v(a - b)$ for all $n \in \Z \setminus \{ 0 \}$.
If there exist $c, d \in D$ such that  $c^n = a, d^n = b$, $\overline{c} = \overline{d}$  for some $n \in \N$, then 
$v( c^m - d^m  ) = v(a - b)$ for all $m \in \Z$. 
\end{lema}
\proof
For $n \in \N$, write 
$a^n - b^n = (a - b)\sum_{i = 0}^{n-1}a^{n- 1 - i}b^i + \text{terms with higher $v$-value}$. Since $\overline{ a^{n- 1 - i}b^i  }$ is the same for all $ 0 \le i \le n - 1$, the $v$-value of the sum is equal to zero, proving the statement for positive integers $n$. For negative $n \in \Z$, the statement follows from
$a^n - b^n = - a^{n}( a^{-n} - b^{- n }) b^n$. The last statement of the lemma follows from
$a - b \sim (c - d    ) \sum_{i = 0}^{n-1} c^{n - 1 - i}  b^i.$
\endproof

\begin{lema} \label{babyconstruct}
Suppose $v$ is a valuation on $\rweyl$ and suppose $i_1, i_2, \ldots, i_r \in \N $ and $k_0 \in \Z$, $k_{i_1 }, k_{i_2 }, \ldots, k_{i_r } \in \Z \setminus \{ 0 \}$ are such that 
$v( x^{k_0}\omega_{i_1 - 1}^{k_{i_1 }} \cdots \omega_{i_r - 1}^{k_{i_r }}    ) = 0$. 
 If $\min_{1 \le j \le r} \{ v(\omega_{i_j }) \}$ is achieved at exactly one $j$, then 
$$v(  x^{k_0} \omega_{i_1 - 1}^{k_{i_1 }} \cdots \omega_{i_r - 1}^{k_{i_r }}     - 
 \overline{ x^{k_0} \omega_{i_1 - 1}^{k_{i_1 }} \cdots \omega_{i_r - 1}^{k_{i_r }}  }  ) 
 = \min \{ v(\omega_{i_j })  \mid 1 \le j \le r \}. $$
\end{lema}
\proof
Let $n$ be the least common multiple of $n_{i_1}, \ldots, n_{i_r}$ and $c_{i_j} = \frac{n}{n_{i_j}}$ for each $i_j$.
Since
\begin{align*}
\overline{(x^{k_0}\omega_{i_1 - 1}^{k_{i_1 }} \cdots \omega_{i_r - 1}^{k_{i_r }})^n} = 
\overline{\Pi_{j= 1}^r x^ { n k_{i_j }  \frac{m_{i_j }}{n_{i_j }}   }(\omega_{i_1 - 1}^{k_{i_1 }} \cdots \omega_{i_r - 1}^{k_{i_r }})^n} = \overline{  \Pi_{j = 1 }^r  ( x^{m_{i_j}} \omega_{i_j - 1}^{ n_{i_j}  }      )^{ k_{i_j}c_{i_j}  }   }
= \Pi_{j = 1}^r \beta_{i_j}^{ k_{i_j}c_{i_j} },
\end{align*}
by Proposition \ref{ratdepelt}, we can compute
\begin{align*}
&(x^{k_0}\omega_{i_1 - 1}^{k_{i_1 }} \cdots \omega_{i_r - 1}^{k_{i_r }})^n     -   \Pi_{j = 1}^r \beta_{i_j}^{ k_{i_j}c_{i_j} }
\\
& \sim \sum_{j = 1}^r  
 ( x^{m_{i_1 } }\omega_{i_1 - 1}^{n_{i_1 }})^{{k_{i_1 } c_{i_1 }}  } \cdots  
 ( x^{m_{i_{j - 1} } }\omega_{i_{j - 1} - 1}^{n_{i_{j - 1} }})^{{k_{i_{j - 1} } c_{i_{j - 1} }}  } \cdot   \\
 &(  (x^{m_{i_j }} \omega_{i_j - 1}^{n_{i_j }} )^{{k_{i_j } c_{i_j }}  }    - \beta_{i_j }^{{k_{i_j } c_{i_j }}  })\beta_{i_{j + 1} }^{  k_{i_j}c_{i_j}  }  \cdots \beta_{i_r }^{ k_{i_j}c_{i_j} }.
\end{align*}
For each $j$ such that $k_{i_j} > 0$,
\begin{align*}
(x^{m_{i_j }} \omega_{i_j - 1}^{n_{i_j }} )^{{k_{i_j } c_{i_j }}  }    - \beta_{i_j }^{{k_{i_j } c_{i_j }}  } & \sim \omega_{ i_{j }} \sum_{i = 0}^{{k_{i_j } c_{i_j }} - 1  } (x^{m_{i_j }} \omega_{i_j - 1}^{n_{i_j }} )^{{k_{i_j } c_{i_j }}  - i }\beta_{i_j }^i,
\end{align*}
which gives us $v( (x^{m_{i_j }} \omega_{i_j - 1}^{n_{i_j }} )^{{k_{i_j } c_{i_j }}  }    - \beta_{i_j }^{{k_{i_j } c_{i_j }}  }   ) = v(\omega_{i_j})$. In case $ k_{i_j} < 0$, we see that
\begin{align*}
   (x^{m_{i_j }} \omega_{i_j - 1}^{n_{i_j }} )^{{k_{i_j } c_{i_j }}  }    - \beta_{i_j }^{{k_{i_j } c_{i_j }}  } &= 
-  (x^{m_{i_j }} \omega_{i_j - 1}^{n_{i_j }} )^{{k_{i_j } c_{i_j }}  }   \beta_{i_j }^{{k_{i_j } c_{i_j }}  }    ((x^{m_{i_j }} \omega_{i_j - 1}^{n_{i_j }} )^{{ - k_{i_j } c_{i_j }}  }    - \beta_{i_j }^{{- k_{i_j } c_{i_j }}  }),
\end{align*}
which again implies $v( (x^{m_{i_j }} \omega_{i_j - 1}^{n_{i_j }} )^{{k_{i_j } c_{i_j }}  }    - \beta_{i_j }^{{k_{i_j } c_{i_j }}  }   ) = v(\omega_{i_j})$.
We conclude, using Lemma \ref{potence}, 
\begin{align*}
v(x^{k_0}\omega_{i_1 - 1}^{k_{i_1 }} \cdots \omega_{i_r - 1}^{k_{i_r }}         -    \overline{x^{k_0}\omega_{i_1 - 1}^{k_{i_1 }} \cdots \omega_{i_r - 1}^{k_{i_r }}}     ) &= 
v( (x^{k_0}\omega_{i_1 - 1}^{k_{i_1 }} \cdots \omega_{i_r - 1}^{k_{i_r }})^n - \overline{(x^{k_0}\omega_{i_1 - 1}^{k_{i_1 }} \cdots \omega_{i_r - 1}^{k_{i_r }}})^n) \\ &= 
 \min \{ v(\omega_{i_j })  \mid 1 \le j \le r  \}  .
\end{align*}
\endproof
With the help of Lemma \ref{babyconstruct}, we will evaluate $v( x^{k_0}  \omega_{i_1 - 1}^{k_{i_1 }} \cdots \omega_{i_r - 1}^{k_{i_r}}     -  \overline{ x^{k_0}\omega_{i_1 - 1}^{k_{i_1 }} \cdots \omega_{i_r - 1}^{k_{i_r}}   }      ) $ 
when $v( x^{k_0}\omega_{i_1 -1}^{k_{i_1 }} \cdots \omega_{i_r -1}^{k_{i_r }}    ) = 0 $ in general.  As in Lemma
\ref{babyconstruct}, we assume $k_0 \in \Z$, $k_{i_j} \in \Z \setminus \{ 0 \}$ for all $1 \le j \le r$. 
This will be helpful when we will later construct a valuation $v$ associated to a 
sequence $( \omega_i)_{i \ge - 1}$. Let us assume for now that $i_1 < i_2 < \cdots < i_r$; at the end of the calculation we will see that the order of $i_j$ does not affect the $v$-value.

To start, we introduce some abbreviations to make the written equations easier to read. Let $n$ and $c_{i_j }$ for all $j$ be as in the proof of Lemma \ref{babyconstruct},
\begin{align*}
A_0 &:=  x^{k_0}\omega_{i_1 - 1}^{k_{i_1 }} \cdots \omega_{i_r -1}^{k_{i_r }}     \\
B_0 &:= \overline{A_0 } \\
A_0^{(n)} &:=x^{nk_0}\omega_{i_1 - 1}^{nk_{i_1 }} \cdots \omega_{i_r - 1}^{nk_{i_r }} \\
B_0^{(n)} &:= \overline{A_0^{(n)}} = \Pi_{j = 1}^r \beta_{i_j }^{ k_{i_j } c_{i_j } }.
\end{align*}
Since $B_0$ is in  $\R$, we can write
\begin{align}\label{enacba1}
A_0 - B_0 &=  (A_0^{n} - B_0^{n}) (\sum_{i = 0}^{n - 1}A_0^{n - i - 1}B_0^i)^{-1} \sim   (A_0^{(n)} - B_0^{(n)}) (\sum_{i = 0}^{n - 1}A_0^{n - i - 1}B_0^i)^{-1}.
\end{align}
Since $v(\sum_{i = 0}^{n - 1}A_0^{n - i - 1}B_0^i) = v(A_0^{n - i - 1}B_0^i) = 0$ for all $i$, which holds due to $\overline{A_0^{n - i - 1}B_0^i} = B_0^{n-1}$ for all $i$, $v(A_0 - B_0) = v(A_0^{(n)} - B_0^{(n)})$. To evaluate the right-hand side of \eqref{enacba1}, we first proceed as we have done in the proof of Lemma \ref{babyconstruct}, so 
\begin{align*}
A_0^{(n)} - B_0^{(n)} &\sim \sum_{j = 1}^r   \Pi_{\ell = 1}^{j - 1} ( x^{m_{i_{\ell} } }\omega_{i_{\ell} -1 }^{n_{i_{\ell}  }})^{{k_{i_{\ell} } c_{i_{\ell}  }}  }  \cdot  (  (x^{m_{i_j }} \omega_{i_j - 1}^{n_{i_j }} )^{{k_{i_j } c_{i_j }}  }    - \beta_{i_j }^{{k_{i_j } c_{i_j }}  }) \cdot
  \Pi_{\ell = j + 1}^r    \beta_{i_{\ell}}^{ k_{i_j } c_{i_j } }. 
\end{align*} If $ k_{i_j} > 0$, we proceed by
\begin{align*}
  (x^{m_{i_j }} \omega_{i_j - 1}^{n_{i_j }} )^{{k_{i_j } c_{i_j }}  }    - \beta_{i_j }^{{k_{i_j } c_{i_j }}  } \sim
\omega_{i_j  } 
\sum_{p = 0}^{{k_{i_j } c_{i_j }} - 1  }  (x^{m_{i_j }} \omega_{i_j - 1}^{n_{i_j }} )^{{k_{i_j } c_{i_j }}  - p - 1 }\beta_{i_j }^p.
\end{align*}
For $i_j > 0, k_{i_j} < 0$, we can on the other hand write 
\begin{align*}
(x^{m_{i_j }} \omega_{i_j - 1}^{n_{i_j }} )^{{k_{i_j } c_{i_j }}  }    - \beta_{i_j }^{{k_{i_j } c_{i_j }}  } \sim 
- \omega_{i_j} (x^{m_{i_j }} \omega_{i_j - 1}^{n_{i_j }} )^{{k_{i_j } c_{i_j }}  }\beta_{i_j }^{{k_{i_j } c_{i_j }}  }
\sum_{p = 0}^{{-k_{i_j } c_{i_j }} - 1  }  (x^{m_{i_j }} \omega_{i_j - 1}^{n_{i_j }} )^{{-k_{i_j } c_{i_j }}  - p - 1 }\beta_{i_j }^p.
\end{align*}
We now define, if $k_j > 0$,
\begin{align*}
C_{ j } &:= \sum_{p = 0}^{{k_{i_j } c_{i_j }} - 1  } 
\Pi_{{\ell} = 1}^{j - 1} (x^{m_{i_{\ell} }} \omega_{i_{\ell} - 1}^{n_{i_{\ell} }} )^{{k_{i_{\ell} } c_{i_{\ell} }}  }   \cdot
   (x^{m_{i_j }} \omega_{i_j -1}^{n_{i_j }} )^{{k_{i_j} c_{i_j }}  - p - 1}\beta_{i_j}^p   
 \cdot  \Pi_{{\ell} = j + 1}^r \beta_{i_{\ell} }^{ k_{i_{\ell} } c_{i_{\ell} } },
\end{align*}
and, if $k_j < 0$,
\begin{align*}
C_{ j } &:= -  (x^{m_{i_j }} \omega_{i_j - 1}^{n_{i_j }} )^{{k_{i_j } c_{i_j }}  }\beta_{i_j }^{{k_{i_j } c_{i_j }}  }  \sum_{p = 0}^{-{k_{i_j } c_{i_j }} - 1  } 
\Pi_{{\ell} = 1}^{j - 1} (x^{m_{i_{\ell} }} \omega_{i_{\ell} - 1}^{n_{i_{\ell} }} )^{{k_{i_{\ell} } c_{i_{\ell} }}  }   \cdot
   (x^{m_{i_j }} \omega_{i_j -1}^{n_{i_j }} )^{-{k_{i_j} c_{i_j }}  - p - 1}\beta_{i_j}^p   
 \cdot  \Pi_{{\ell} = j + 1}^r \beta_{i_{\ell} }^{ k_{i_{\ell} } c_{i_{\ell} } }
\end{align*}
for each $1 \le j \le r$. Further, we define
\begin{align*}
A_{1, j } &= \omega_{i_j }  C_{j} \\    
A_0^{(n)} - B_0^{(n)} &\sim \sum_{j = 1}^r  \omega_{i_j } C_{j}   =   \sum_{j = 1}^r A_{1, j}
\end{align*}
for each $1 \le j \le r$. Thus we can conclude that $v(A_{1, j}) = v( \omega_{i_j } )$, since the image of $C_{ j}$  in the residue field is equal 
$  k_{i_j  } c_{i_j }  \cdot \Pi_{{\ell} = 1}^r \beta_{i_{\ell} }^{k_{i_{\ell}} c_{i_{\ell} } } \neq 0$ if $k_{i_j} > 0$ and
$ - \beta_{i_j}^{2 k_{i_j}c_{i_j}} |k_{i_j  } c_{i_j }|  \cdot \Pi_{{\ell} = 1}^r \beta_{i_{\ell} }^{k_{i_{\ell}} c_{i_{\ell} } } \neq 0$ if $k_{i_j} < 0$
, making $v(C_{j})  = 0$. We can now write
\begin{align*}
A_0^{(n)} - B_0^{(n)} &= \sum_{j = 1}^r A_{1, j} + A = \sum_{v(A_{1, j}) \text{ is minimal} } A_{1,j} + \sum_{v(A_{1, j}) \text{ is not minimal} } A_{1, j}.
\end{align*}
Here we note that the second of both finite sums on the right-hand side of this equation includes $A$, which denotes the sum of all terms obtained by changing the order of factors of the form $\omega_{i_{\ell }}  $ (which was not explicitly written above). The fact that the $v$-value of these terms is higher than the $v$-value of the terms of the first sum (the ones with minimal $v$-value)  follows from Lemma \ref{kompogoj}.

If $\min_{1 \le j \le r} \{ v(\omega_{i_j }  ) \} $ is achieved at more than one $j$, we take the sum of all $ \omega_{i_j } { C_{j}   }$ that have the minimal $v$-value, i.e., $\sum_{v(A_{1, j}) \text{ is minimal} } A_{1,j}$, then factor $ \omega_{i_1 } $, so the sum now looks like 
$$  \sum_{v(A_{1, j}) \text{ is minimal} } A_{1,j} =   \omega_{i_1 } \sum_{v(A_{1, j}) \text{ is minimal} }\omega_{i_1 }^{-1}  A_{1,j}  =      \omega_{i_1 }  \sum_{v(A_{1, j}) \text{ is minimal} } \omega_{i_1 }^{-1} \omega_{i_j } { C_{j}   }$$
 and, since $v (\omega_{i_1 }^{-1} \omega_{i_j }C_j) = 0$, for each $j$, we can evaluate the sum of their images in the residue field. If this sum is not equal to zero, then
$v(A_0 - B_0)  = \min_{1 \le j \le r} \{ v(A_{1,j}) \} $. 
Otherwise write
\begin{align*}
\omega_{i_1 } \sum_j  \omega_{i_1 }^{-1} \omega_{i_j } { C_{j}   } =
\sum_j \omega_{i_1 } ( \omega_{i_1 }^{-1} \omega_{i_j } { C_{j}   } - \overline{  \omega_{i_1 }^{-1} \omega_{i_j } { C_{j}   }   }).
\end{align*}
For every $j$, we write $  \omega_{i_1 }^{-1} \omega_{i_j } { C_{j}   } - \overline{  \omega_{i_1 }^{-1} \omega_{i_j } { C_{j}   }   }$ as an $\R$-linear sum of terms of the form $\Pi_{\ell} \omega_{i_{\ell} - 1 }^{n_{\ell}}$ (in the same way we did with $A_0 - B_0$). We sum all of the newly obtained terms, as well as the terms in $ \sum_{v(A_{1, j}) \text{ is not minimal} } A_{1, j}   $, and relabel them as $A_{2, j}$ where $j$ goes from $1$ to the number of all terms. 

As $A_0 - B_0$ can be written in the form 
\begin{align*}
A_0 - B_0 = (\sum_{v(A_{2,j}) \text{ is minimal }   }  A_{2,j} +  \sum_{v(A_{2,j}) \text{ is not minimal }   }  A_{2,j} )D,
\end{align*}
where we use $D$ as the label of the product of all terms  of the form $( \sum_i A^{k - i}B^i)^{-1}$ where
$A = \overline{x^{k'_0} \Pi_{j = 1}^r \omega_{i_j-1}^{k'_{i_j}}  }$ for some $i_1, \ldots, i_r \in \N$, $k'_0, k'_{i_1}, \ldots, k'_{i_r} \in \Z$ and $B = \overline{A}$, that we factor out when we evaluate $  \omega_{i_1 }^{-1} \omega_{i_j } { C_{j}   } - \overline{  \omega_{i_1}^{-1} \omega_{i_j } { C_{j}   }   }$ for each $j$. All terms $( \sum_i A^{k - i}B^i)^{-1}$ have $v$-value equal to zero and their image in the residue field, which is of the form 
$(m B^{m - 1})^{-1} \in \R$ for some $m \in \N$,  is easy to determine.

We repeat the described procedure, writing $A_0 - B_0 = (\sum_j A_{k,j})D$ for increasing $k$. We stop when for some $k$, $\sum_{j, v(A_{k, j}) \text{ is minimal}   }A_{k, j}$ is either composed of one single term or, after factoring out one of the terms, the image of the sum in the residue field is not zero. In this case, we conclude that $v(A_0 - B_0)$ is equal to $v( A_{k, j} )$ for any term of the sum $\sum_{j, v(A_{k, j}) \text{ is minimal}   }A_{k, j}$.

We must show that the process ends at some point even if the number of terms whose $v$-value we evaluate at each step is growing. We see that whenever we write $x^{k_0}\Pi_{\ell = 1}^r \omega_{i_{\ell} - 1}^{k_{\ell} } - \overline{x^{k_0}\Pi_{\ell = 1}^r \omega_{i_{\ell} - 1}^{k_{\ell } }   }   $ as a sum of terms with strictly positive $v$-value, the value of each of these terms is $v(\omega_{i_{\ell}  })$ for some $\ell = 1, \ldots, r$. It follows that $v(A_0 - \overline{A_0})$ is a sum of $v( \omega_{\ell}  )$ for some $\ell \ge 1$.

If $v(\omega_N)$ is irrational for some $N \ge 0$, the process either stops beforehand or, after $k \ge N - i_r$ steps we get a unique term $A_{k, j}$ that has $v$-value equal to $v(\omega_{i_r } \omega_{i_r + 1} \cdots \omega_N)$. 
This is the term we get when we take the last term of $A_0^{(n)} - B_0^{(n)}$, written as a sum of terms $A_{1, j}$ with higher $v$-value and in each of the following steps whenever the $v$-value of this term is minimal, take the last term when $A_{k,j}$ is written as a sum of terms with higher $v$-value. 

If on the other hand, $v(\omega_k) \in \Q$ for an infinite sequence $( \omega_k  )_{k \ge -1}$, then
 $\lim_{k \rightarrow \infty} v( \omega_k)  = 0$ since by Lemma \ref{kompogoj}, 
$\sum_{i \ge -1}^k v( \omega_k )< 0$ for all $k \ge -1$ and $v(  \omega_i  ) > 0   $ for $i \ge 1$.
Then for some $N \ge 1$, $v(    \omega_N ) < v( \omega_i  )$ for all $1 \le i < N$. The evaluation of 
$v(A_0 - B_0)$ again either stops beforehand or we get a unique term $A_{k, j}$ that has $v$-value equal to 
$v(\omega_{i_r } \omega_{i_r +1 } \cdots \omega_N)$. As in the first case, this term is the one we get when we take the last term of $A_0^{(n)} - B_0^{(n)}$, written as a sum of terms $A_{1, j}$ with higher $v$-value and in each of the following steps whenever the $v$-value of this term is minimal, take the last term when $A_{k,j}$ is written as a sum of terms with higher $v$-value. 

In both cases, the value of this term, $v( \omega_{i_r } \omega_{i_r + 1} \cdots \omega_N  )$  is strictly smaller than the value of all additional terms we get when we change the order of the factors in a product. It follows that given the $v$-values of $(\omega_i)_{i \ge -1}$, $v(A_0 - B_0)$ is the same as it would be if all elements of the sequence $(\omega_i)_{i \ge -1}$ commuted.
This follows from the fact that when we change the order of factors $\omega_i$ and $\omega_j$ in some $A_{k, j_1}$, the term we obtain has $v$-value greater by $v[\omega_i, \omega_j ] - v( \omega_i \omega_j) = - v(xy \omega_1 \cdots \omega_j)$ for $-1 \le i < j \le N$ while $v(A_{k + 1, j_1}) - v(A_{k, j_2})$ is for any $j_1, j_2$ equal to $v(\omega_{i_{\ell} +  1 })$ for some $i_{\ell}$ such that $A_{k, j _2}$ contains a power of $\omega_{i_{\ell}}$. 
Since, as we have shown in the proof of Lemma \ref{kompogoj}, $\sum_{i = -1}^N v(\omega_i) < 0$, it follows that the terms we obtain by changing the order of the factors have $v$-value greater than $v(A_0 - B_0)$. 
And since $v(   x^{k_0}\omega_{i_1 - 1}^{k_{i_1 }} \cdots \omega_{i_r - 1}^{k_{i_r}}     - 
\overline{x^{k_0}\omega_{i_1 - 1}^{k_{i_1 }} \cdots \omega_{i_r - 1}^{k_{i_r}}}   )$ is higher than the $v$-value of any terms we get when we change the order of factors, the order of $\omega_{i_1}, \ldots, \omega_{i_r}$ does not matter. 

\subsection{Extending $v$ from the sequence $ ( \omega_i)_{i \ge -1} $ to $\rweyl$ }

We can now prove that every $v$ associated to either a finite or an infinite sequence $( \omega_i)_{i \ge - 1}$ can be extended to a valuation on $\mathcal{D}_1( \R)$.

\begin{lema} \label{diophantine}
For every $r \ge 0$, there exists a finite number $M$ of elements of the form $a_i := \omega_{ - 1}^{k_{i, -1}} \omega_0^{k_{i,0}} \omega_1^{k_{i, 1}} \cdots \omega_{r - 1}^{k_{i, r  - 1}}$ for some $k_{i, 1}, \ldots, k_{i, r- 1} \in \Z$  such that $v(a_i) = 0$ for all $1 \le i \le M$ and every
$ \omega_{-1}^{\ell_{-1}} \cdots \omega_{r-1}^{ \ell_{r - 1}  }  \in \mathcal{D}_1( \R )$ with $v$-value zero is $\sim$-equivalent to a product of positive integer powers of $a_1, \ldots, a_M$.
\end{lema}

\proof
Since $v(\omega_{-1}) = -1$ and    $v(\omega_{j}) = \frac{m_{j + 1}}{ n_{j + 1}   } \in \Q $ for $j \ge 0$, the problem translates to finding general classes of solutions to 
the diophantine equation $x_0 a_0 + \cdots x_k a_k = 0$ with $a_0 = - \Pi_{j = 1}^k n_i$ and $a_i = \frac{m_i}{n_i} \Pi_{j = 1}^k n_i$ for all $0 \le i \le k$. 
\endproof

\begin{izrek} \label{konstrukcija2}
Let $v$ and $(\omega_i)_{ i \ge -1 }$ be as described in the beginning of the section, i.e., $\omega_{-1} = x$, $\omega_0 = y$, $v( \omega_{-1}) = -1$, $v(\omega_i) = \frac{m_{i + 1}}{ n_{i + 1}   }  \in \Q$, $\overline{x^{m_{i + 1}} \omega_i^{ n_{i + 1} } } = \beta_{i + 1} \in \R  $, 
$\omega_{i +1} = x^{m_{i + 1}} \omega_i^{ n_{i + 1} } - \beta_{i + 1} $ for $i,  0 \le i \le N - 1$ and 
$v(\omega_N) \not\in \Q$ for some $N \ge 0$ or $v(\omega_i) \in \Q$ for infinitely many $i$, that $\sum_{i = -1}^k v( \omega_i  ) < 0$ for all $k \ge - 1$.
Suposse that $\sgn{\beta_i}$ is constant on the set of all $i$ for which $n_i$ is even. 
Then $v$ can be extended to a valuation on $\mathcal{D}_1( \R )$ with residue field $\R$. The valuation is unique for every choice of $\{  \alpha_{i, j}\}_{i, j \ge 0}$ where 
$\alpha_{i, j} = \overline{\omega_{i - 1}^{K_{i,j} }\omega_{j - 1}^{K_{j, i}}}, K_{i, j} = \frac{m_jn_i}{d_{i,j}}, K_{j,i} = - \frac{m_in_j}{d_{i,j}}
\text{with }d_{i,j} = \gcd \{ m_jn_i, m_in_j \}.$
The associated value group is group-isomorphic to a subgroup of $\Q \times \Z$ generated by $\{ v(\omega_i)  \}_{i \ge -1}$.
\end{izrek}

\proof
The following construction of the valuation $v$ associated to the sequence $(\omega_i)_{i \ge -1} $  was first sketched in \cite{SHT}. Here we present it in full detail. 

Before we begin with the construction of the $v$-value for an arbitrary element of $\mathcal{D}_1 ( \R )$, we define it for some specific elements of $\mathcal{D}_1 (\R)$.
\begin{enumerate}
\item Since we have defined $v( \omega_i  )$ for all $-1 \le i \le N$, 
$v( \Pi_{i = - 1}^N \omega_i^{k_i}    ) =  \sum_{i = -1}^N k_i v( \omega_i  )    $ must hold for all 
$k_{-1}, \ldots, k_N \in \Z$.
\item Since we supposed $\sum_{i = -1}^k v( \omega_i  ) < 0$ for all $k \ge - 1$, 
it follows from Lemma \ref{kompogoj} that 
$$v[ \omega_i, \omega_j ] = - v(\omega_{-1} \omega_0 \cdots \omega_{i-1}\omega_{i + 1}\cdots \omega_{j - 1}  ),$$ which must be strictly greater than $v(\omega_i \omega_j)$ for all $i < j$. 
\item We also see that if $v( \omega_{-1}^{k_{0}} \cdots \omega_{r - 1}^{k_r}   ) = 0$ for some $k_0, \ldots, k_r \in \Z$, then
$\overline{  \omega_{-1}^{k_{0}} \cdots \omega_{r - 1}^{k_r}   }$ is uniquely determined by $\{ \alpha_{i,j} \}_{i, j \ge 0}$ as in shown in Proposition \ref{dvojicedet}. In this case, $v( \omega_{-1}^{k_{0}} \cdots \omega_{r - 1}^{k_r}       - \overline{ \omega_{-1}^{k_{0}} \cdots \omega_{r - 1}^{k_r}    }     )$  must be equal to the value determined in Lemma \ref{babyconstruct}  and the discussion following it. 
\end{enumerate}
In all three cases, the chosen values were the only possible extensions of $v$ from $( \omega_i)_{i \ge - 1}$ if we want $v$ to be a valuation.

To determine $v(F)$ for any $F \in \rweyl$, we first note that $F$ can be written as a finite sum 
\begin{align*}
F = \sum_{\ell } \alpha_{\ell}x^{i_{\ell}   }y^{j_{\ell}}, \alpha_{\ell} \in \R.
\end{align*}
Let $F_1$ be the sum of all terms $\alpha_{\ell}x^{i_{\ell}   }y^{j_{\ell}}$ such that $i_{\ell} v( x  ) + j_{\ell}v(y) $ is equal to $ u := \min_{ \ell } \{  i_{\ell} v( x  ) + j_{\ell}v(y)  \}.    $

 If $F_1$ consists of only one such term, then we define $v(F) = u$; this is obviously always the case whenever $v(y) \not\in \Q$. Otherwise, we factor out $x^{i_1}y^{j_1}$ with the smallest power of $x$ and get

\begin{align*}
F_1 
 \sim x^{i_1}y^{j_1}  \sum_{ \ell, \  i_{\ell} v( x  ) + j_{\ell}v(y)  = u } \alpha_{\ell}x^{i_{\ell} - i_1  }y^{j_{\ell} - j_1}.
\end{align*}
Since  $$(i_{\ell} - i_1 )v(x)  + (j_{\ell} - j_1)v(y) = 0, $$ for each $\ell$ in the sum, $\frac{i_{\ell} - i_1}{j_{\ell} - j_1} =   K_{\ell}  \frac{ m_1  }{ n_1}  $ for some 
$K_{\ell} \in \Z$. We can
write $$F_1 \sim x^{i_1}y^{j_1} f( x^{m_1}y^{n_1}   ) $$ where $f(t)$ is a polynomial in $\R[t]$. 
Since we know $v(\omega_1) > 0$ and $v(\alpha) = 0$ for $ \alpha \in \R^*$, $v$ is uniquely determined on $\R[ \omega_1]$. From this, it follows that $v(f(x^{m_1}y^{n_1}    )) = 0$ if and only if $f(\beta_1) \neq 0$ since $x^{m_1}y^{n_1} = \omega_1 + \beta_1  $. 

In this case, $v(F_1) = u$ and since all terms in $F - F_1$ have $v$-value strictly greater than $u$, $v(F) = u$ must hold. 
Since 
$v(u)$ is a sum of integer powers of $v(x )$ and $v(y)$, $v(F)$ is in the abelian group, generated by $\{ v(\omega_i) \}_{i \ge - 1}.$ If $u = 0$, 
$\overline{F} = \beta_1^k f( \beta_1 ) \in \R$. 
If on the other hand $f(\beta_1 )= 0$, write
$f(t) = g_1(t)(t - \beta_1)^{k_1}$ with $g_1(\beta_1) \neq 0$ and we have 
$$F_1 \sim x^{i_1}y^{j_1}g_1(x^{m_1}y^{n_1}) \omega_1^{k_1}    .$$ 
We set $v(F_1) = u + k_1 v( \omega_1)$ and add all terms we get from exchanging the order of factors, whose $v$-value can be lower than the newly set $v(F_1)$, although still strictly higher than $u$ due to $v[x, y] > v(xy)$, to $F - F_1$. It is immediate that $v(F_1) $ is in the subgroup of $\Gamma$, generated by $v(x), v(y)$ and $v(\omega_1)$ and that if $v(F_1) = 0$, $\overline{F_1} \in \R$.  It is important to note that in both cases, we consider $F_1$ as a single term. It follows that during our transformation, the number of terms (if we ignore the ones we got when we changed the order of factors in a product) is strictly smaller than before (unless, of course, $F_1$ was just a single term in the beginning and we get $v(F) = v(F_1)$).
 
We now consider the values of the terms in $F - F_1$. If all of them have $v$-value strictly greater than that of $F_1$, we conclude $v(F) = v(F_1)$. Otherwise, we take all terms of 
\begin{align*}
F - F_1  =  \sum_{\ell', \  i_{\ell'} v( x  ) + j_{\ell'}v(y)  > u } \alpha_{\ell'}x^{i_{\ell'}   }y^{j_{\ell'}}, \alpha_{\ell'} \in \R
\end{align*}
for which $i_{\ell'}v(x) + j_{\ell'}v(y) = u' := \min_{\ell' } \{ i_{\ell'}v(x) + j_{\ell'}v(y) \} $ and then as before define
\begin{align*}
F_2 = x^{i_2}y^{j_2}  \sum_{ \ell', \  i_{\ell'} v( x  ) + j_{\ell'}v(y)  = u' } \alpha_{\ell'}x^{i_{\ell'} - i_2  }y^{j_{\ell'} - j_2}.
\end{align*}
As above, we write $F_2 \sim x^{i_2} y^{j_2}  g_2(x^{m_1}y^{n_1} ) ( x - \beta_1  )^{k_2} $, $k_2 \ge 0$, $g_2(\beta_1) \neq 0$ and add all the terms we get when we change the order of factors to $F - F_1 - F_2$. Their $v$-value is strictly greater than $u'$ due to $v[x, y] > v(xy)$.

We continue this process, defining $F_1, F_2, \ldots, F_k$ until all terms in $F - F_1 - \cdots  - F_k$ have $v$-value strictly greater than 
$\min\{ F_1, \ldots, F_k \}$. Note that it is possible that $F_k$ consists of only one term from $F - F_1 - \cdots - F_{k-1}.$

Afterwards, we sum together all those $F_i$ for $1 \le i \le k$ for which $v(F_i) = u_1 := \min \{ v(F_1), \ldots, v(F_k) \}$. 
If the minimum is achieved at exactly one such $F_i$, we set $v(F) = v(F_i)$. This is always the case whenever $ v(\omega_1) \not\in \Q $. 
 Otherwise we can relabel the terms so the minimum is achieved at 
$F_1, \ldots, F_r$ for some $r \le k$. As we have shown, each $F_i $ can be written as 
$F_i =  x^{i_i} y^{j_i} g_i(x^{ m_1} y^{ n_1  }) \omega_1^{k_i}  $. 
We sum the terms together, factor out $x^{i_1}y^{j_1} \omega_1^{k_1} $, the term that has, written as a polynomial in $x$ and $y$, the lowest power of $x$,  and label the new sum  $F_{1, 1}$. 

To evaluate $v(F_{1,1})$, we follow a procedure similar to the one evaluating $v(F_1)$. After factoring
$x^{i_1}y^{j_1} \omega_1^{k_1}$, we are left with
\begin{align*}
F_{1,1} &\sim x^{i_1}y^{j_1} \omega_1^{k_1}( g_1(x^{m_1}y^{n_1}) +  g_2(x^{m_1}y^{n_1}) x^{i_2 - i_1}y^{j_2 - j_1}  \omega_1^{k_2 - k_1} +
 \cdots \\ 
&+ g_r(x^{m_1}y^{n_1}) x^{i_r - i_1}y^{j_r - j_1}  \omega_1^{k_r - k_1}  ) \\
&\sim x^{i_1}y^{j_1} \omega_1^{k_1} \sum_{J = 1}^R \alpha_J  x^{i_J}y^{j_J}\omega_1^{k_J} \text{, with $\alpha_J \in \R,   i_J, j_J, k_j \in \Z, R \ge 1.$}
\end{align*}
Each term in the sum has $v$-value zero.
Let $a_1, a_2, \cdots, a_{\ell}$ be the terms such that each product of the form  $x^{i}y^{j}\omega_1^k$, $i,j,k \in \Z$ that fulfills the condition $v(x^{i}y^{j}\omega_1^k) = 0$ is a product of positive integer powers of some of $a_i$ up to the order of factors $x, y$ and $\omega_1$. The existence of $a_1, \ldots, a_{\ell}$ is assured by Lemma \ref{diophantine}. We can then write
\begin{align}\label{enacba2}
F_{1,1} &\sim x^{i_1}y^{j_1} \omega_1^{k_1} \sum_{J = 1}^R  \gamma_J a_1^{m_{1, J}}a_2^{m_{2, J}}\cdots a_{\ell}^{m_{\ell, J}} \nonumber \\
&\sim x^{i_1}y^{j_1}  \omega_1^{k_1} g(a_1, \cdots, a_{\ell}), g \in \R[t_1, \cdots, t_{\ell}]
\end{align}
with $\gamma_J \in \R$, $m_{I, J} \in \Z$ for all $1 \le I \le \ell$ and  $ 1 \le J \le R$. As before, we add all terms we get when we change the order of multiplication of $x, y$ or $\omega_1$ in a product to $F - F_{1,1}$ since the value of its terms is strictly greater than $u_1$.  Since, as we have determined in the beginning, each term in the sum \eqref{enacba2}  has $v$-value equal to zero and we know what $\overline{a_i} \in \R$ is for each $1 \le i \le \ell$, $v(g(a_1, \ldots, a_{\ell}))$ will have to be greater than or equal to zero, we can define 
$\overline{g(a_1, \ldots, a_{\ell})} = g(\overline{a_1}, \ldots, \overline{a_{\ell}} ).  $ 

If $g(\overline{a_1}, \ldots, \overline{a_{\ell}} ) \neq 0$, then we set $v(g(a_1, \ldots, a_{\ell})) = 0$
and $v(F) = v(F_{1,1})  = i_1v(x) + j_1 v(y) + k_1 v(\omega_1)$. Otherwise write 
\begin{align*}
g(t_1, \ldots, t_{\ell}) = h(t_1 - \overline{a_1}, \ldots, t_{\ell} - \overline{a_{\ell}}     ) = \sum_{i = 1}^L \Pi_{j = 1}^ {\ell} (t_j - \overline{a_j})^{m_{i, j}  } h_i(t_1, \ldots, t_{\ell})
\end{align*}
with $h, h_1, \ldots, h_L \in \R[t_1, \ldots, t_n]$ and  $h_i(\overline{a_1}, \ldots, \overline{a_{\ell}} ) \neq 0$ for all $i$. We factor out the $\Pi_j (t_j - \overline{a_j})^{m_{i,j}}$ for those $i$ for which $\sum_j v(a_j - \overline{a_j})^{m_{i,j}}$ is minimal. Then 
\begin{align*}
g(t_1, \ldots, t_{\ell}) = \Pi_j (t_j - \overline{a_j})^{m_{i,j}} \tilde{g}(t_1, \ldots, t_k), \text{ with $\tilde{g} \in \R[t_1, \ldots, t_{\ell}]. $ }
\end{align*}
If $\tilde{g}(\overline{a_1}, \ldots, \overline{a_k}) \neq 0$, we set 
$$v(F_{1,1}) = v( x^{i_1}y^{j_1}  ) + \sum_{j} m_{i,j }v(a_j - \overline{a_j}   ). $$
If on the other hand, $\tilde{g}(\overline{a_1}, \ldots, \overline{a_k}) = 0$, we do the same thing as we did with $g$. The process cannot go on indefinetly since $g$ is a polynomial and hence of finite degree. 
All terms we get when we exchange the order of $x, y$ and $\omega_1$ are added to $F - F_{1,1}$. Their $v$-value must be  strictly greater than $u_1$. It follows from the construction that $v(F_{1,1})$ must be in the group generated by $ \{ v(\omega_i) \}_{i \ge -1} $ since this holds for $v(a_i - \overline{a_i}    )$ for all $i$ and that if $v( F_{1,1}) = 0$, $\overline{F_{1,1}} \in \R.$

Since $v( a_i - \overline{a_i}  )$ is, as we have shown in Lemma \ref{babyconstruct} and the discussion following it, a sum of $v(\omega_j)$ and thus $v(\Pi_j \omega_j^{-1} (a_i - \overline{a_i} )  ) = 0$, we can write $F_{1,1}$ as one term of the form $\Pi_{i = -1}^n \omega_i^{k_i} g(a_1, a_2, \ldots, a_{\ell}  )$ with $n \in \N$ and $g \in \R[t_1, \ldots, t_{\ell}]$ and  $ g(\overline{a_1}, \overline{a_2}, \ldots, \overline{a_{\ell}}  )  \neq 0$.

After $v(F_{1,1})$ is set, we compare it to both $v(F_i)$ for all $F_i$ that are not part of $F_{1,1}$ and the terms of $F - F_1 - \cdots - F_k - F_{1,1}$. If 
all of these terms have $v$-value strictly greater than $v(F_{1,1})$, then we can set $v(F) = v(F_{1,1})$. Otherwise, we collect all terms with minimal $v$-value in a sum which we label $F_{1,2}$. We determine $v(F_{1,2})$ in the same way we determined $F_{1,1}$  and then sum all of the remaining terms that have $v$-value less or equal to $\min \{  v(F_{1,1}), v(F_{1,2})     \}$ to a sum labeled $F_{1,3}$.

We repeat the process until for some $k$, $\min \{ v( F_{1,1}), \cdots, v(F_{1,k})    \}$ is strictly smaller than the $v$-value of any of the remaining terms. 

If $\min \{ v( F_{1,1}), \cdots, v(F_{1,k})    \}$  is achieved at exactly one $i$, we set $v(F) = v( F_{1,i})$. Otherwise we sum all the terms with the minimal $v$-value and label the sum $F_{2,1}$. We evaluate $v(F_{2,1})$ in the same way we evaluated $v(F_{1,1})$. 
 We repeat the process, defining $ F_{i,j} $ and determining its $v$-value in the same way as above. We point out that after $v(F_{i, j})$ is defined, we regard  $F_{i, j }$ as one single term in future evaluations.

Now we must show that at one point, the process ends, i.e., that for some $i, j$, $v(F_{i,j})$ is strictly smaller than the $v$-value of all other terms. This holds because each time we define $F_{i,j}$ for some $i,j$, we sum a number of different terms into one single term and because whenever we change the order of factors in a term, the degree of $x$ and $y$ in the difference is strictly smaller. This means that we eventually run out of terms.
We have thus defined $v$ for an arbitrary polynomial $F \in \rweyl$. What we essentially did was that we wrote
$$F = \tilde{F} + \tilde{F_1} $$
where $\tilde{F}$ is written as a single term, $v(\tilde{F})$ is computed as if $x$ and $y$ commuted and the $v$-value of each term of $\tilde{F_1} $ is strictly greater than $v(\tilde{F})$. For another $G \in \rweyl$, we can write
$$FG = \tilde{FG} + \tilde{(FG)_1}$$
and since we evaluate $v(\tilde{F})$ and $v( \tilde{G})$ as if $x$ and $y$ commuted, $v(\tilde{FG}) = v(\tilde{F}  ) + v(\tilde{G})$. We use the same reasoning to show $v(F + G) \ge \min \{v(F), v(G) \}$.

It follows from the construction that for each $i, j$, $v(F_{i, j })$ is a linear combination of $\{ v( \omega_i  ) \}_{i \ge -1}$ and that in case $v( F_{i, j}   ) = 0$, $\overline{F_{i,j}} \in \R .$ 
\endproof

\begin{izrek} \label{stronglycommutative}
Let $v$ be a valuation on $\rweyl$ trivial on $\R$ with residue field $\R$. Then $v$ is strongly abelian.
\end{izrek}
\proof
If $v$'s value group is $\Q$, then the theorem follows from Corollary \ref{ratvalgrp}. Otherwise, $v(\omega_N) \not \in \Q$ for some $N$ by our construction. 
But as we have shown in Lemma \ref{kompogoj}, $v[\omega_N, x  ] = - v(y \omega_1 \cdots \omega_{N - 1})$. If $v(\omega_N x) = v[ \omega_N, x]$, it follows that $v(\omega_N) \in \Q$, a contradiction. Since the value group is generated by $ \{ v(\omega_i) \}_{i \ge -1}$, it follows from Proposition  \ref{realrank2} that $v$ is strongly abelian.
\endproof

\section{Valuations on $R[y ; \delta]$} \label{middle}

In this section, we explain a construction of valuations on the ring $R[y ; \delta ]$ with $$R := \{ \sum_{k \ge m} a_k x^{- \frac{k}{n} }  \mid a_k \in \R, m \in \Z, n \in \N  \}$$ and  $\delta(p(x)) = p'(x)$. This construction, which was first introduced in \cite{MA2}, will, as we will see in this section, give us all valuations on $R[y ; \delta ]$ with residue field $\R$. Then, we will prove exactly which valuations on $\rweyl$ with residue field $\R$ extend to a valuation on $R[y ; \delta ]$ with the same residue field, answering the question posed by Marshall and Zhang in \cite{MA2}.
We will see the extensions of valuations on $R[y ; \delta]$ are strongly abelian.

Every valuation on $R[y ; \delta]$ can be uniquely extended to its quotient ring, which we label as $D$, because $ R[y ; \delta]$ is an Ore domain.
 Since $[y, x ] = 1$ as before, $v(xy) < 0$ must hold. We set $v(x) = - 1$, $z_0 := y$ and consider $v(y)$. If $v(y) \not \in \Q$, then
$$v( \sum_{i = 0}^n p_i(x) y^i   ) =  \min_{0 \le i \le n} \{ v(p_i(x)) + i v(y)   \}  $$ for any 
$\sum_{i = 0}^n p_i(x) y^i  \in R[y ; \delta]$. Otherwise $v(y) = r_1 \in \Q$ and hence $v(y -  \gamma_1 x^{-r_1}) > v(y)$ for some 
$\gamma_1 \in \R$.  If $v( z_1) = r_2 \in \Q$ for $z_1 := y - \gamma_1 x^{- r_1} $, we proceed to find $z_2 = z_1 - \gamma_2x^{- r_2}$ such that $ v(z_2)$ is greater than $r_2$. We repeat this process to construct a sequence $( z_i )_{i \ge 0}$.

If $v( z_k ) \not\in \Q$ for some $k \in \N$, then we can write every $f \in R[y ; \delta]$ as $ \sum_{i = 0}^n p_i(x)   z_k^i   $ and deduce
$$ v(f ) =   \min_{0 \le i \le n} \{ v(p_i(x)) + i v (z_k)   \}. $$ The value group is then group-isomorphic to $\Q \times \Z$. Since $v[x, z_k ] = v[x, y] = 0 > v(x z_k)$, $v$ is stongly abelian by Proposition \ref{realrank2}.  Otherwise, the sequence 
$( z_i )_{i \ge 0}$, $v( z_i) = r_{i + 1} \in \Q$ is infinite. We take note of the fact that  $v( z_{i + 1}) > v( z_i)$ and since $[ z_i, x] = [y, x] = 1$ for all $i$, $v( z_i) < 1$  for all $i$.
We define $r := \lim_{i \rightarrow \infty}r_i \le 1$. 

\subsection{Case $r < 1$} If $r < 1$, it has been shown in \cite{MA2} that $v$ can be extended to a valuation on ${R}[y ; \delta]$ with residue field $\R$. We first extend $v$ from $R$ to $$\tilde{R} = \{  \sum_{q \in A}a_q x^{-q}  \mid  a_q \in \R,  A \subset \Q \text{  is 
well-ordered}   \}   $$ in a natural way, i.e., by defining 
$$v(  \sum_{q \in A}a_q x^{-q})  = \min A $$
for each $\sum_{q \in A}a_q x^{-q} \in \tilde{R}[y ; \delta]$.
Then for every $f(t) \in R[t]$, define $v(f(y)) = v(f(z))$ with $z := \sum_{i \ge 1}^{\infty} a_i x^{- r_i}$ and $f(y) = \sum_{i = 0}^n p_i y^i $ for $f(t) = \sum_{i = 0}^n p_i t^i $. This gives rise to a valuation on $ R[ y ; \delta ]  $.

However if $r =1$, we cannot define a valuation in this way.
Let $k \in \N$ be such that $ 2 r_k > 1 + r_1  $, which exists since $r =1$, and $a_k = y - z_k = \sum_{i = 1}^k \gamma_i x^{- r_i}$. Let
$$f(t) = (t - a_k)(t - a_k) = t^2 - 2 ta_k + a_k^2 \in R[t].$$
On one hand, $$v(f(y)) = v(f(z)) = 2v(z - a_k) =2 r_{k + 1}.$$ On the other, 
$$2 r_{k + 1} = v ( (y - a_k)(y - a_k)) = v(y^2 - 2ya_k + a_k^2 + [y, a_k]) = \min \{ 2 r_{ k + 1}, 1 + r_1  \} = 1 + r_1 ,$$ contradicting the assumption that $v$ is a valuation, as shown in \cite{MA2}. 

Of course, even in case $r < 1$, there may also exist a $k$ such that $2 r_k > 1 + r_1$. But the important difference between the two cases is that if $r < 1$, there is always an $\ell \in \N$ such that $1 + r_{\ell} > 2 r_k$ for all $k \in \N$, which does not hold in case $r = 1$. Then, since
$$f : R[y ; \delta ] \rightarrow R[y ; \delta], \\ f(y) = z_{\ell}, f(a) = a \text{ for } a \in R$$ is a real algebra automorphism of $R[y ; \delta]$, we can translate the sequence by replacing $y$ with $ z_{\ell}$.

We see that since the associated value group is $\Q$, $v$ is a strongly abelian valuation by Corollary \ref{ratvalgrp}.

\subsection{Case $r = 1$} The question whether in case $r = 1$, $v$ can be extended from a sequence $( z_i )_{i \ge 0}$ to a valuation on  $R[y ;  \delta]$ was left open in \cite{MA2}. In this subsection, we show that it can be done using model theory (for reference, see for example \cite{PRE}). We also show that the valuation we get in this way is uniquely determined.

Suppose we have infinite sequences $( z_i)_{i \ge 0} \subseteq R[y ; \delta]$, $( r_i)_{i \ge 1}  \subseteq \Q$ and $( \gamma_i)_{i \ge 1} \subseteq \R $ and $v : ( z_i)_{i \ge  0} \rightarrow \Q$ with $z_0 =y$, $z_{i + 1} = z_i - \gamma_{i + 1}x^{- r_{i + 1}} $ and $v(z_i) = r_{i + 1} \in \Q$ with $(r_i)_{i \ge 1}$ a strictly increasing sequence with $r = \lim_{i \rightarrow \infty} r_i = 1$.
Then for each $n \ge 0$, there is a valuation $v_n$ on $R[y ; \delta]$ such that $v_n(z_i) = r_{i + 1}$ for all $0 \le i \le n- 1$ and 
$v_n(z_n) \not \in \Q$.

We now present the first-order theory that the valuation associated to the infinite sequence we wish to prove exists is a model of. The theory will be a union of the theory of $D$, the quotient division ring of $R[y ; \delta]$ and the theory of valuations. We will see that each finite subset of this theory has a model. By compactness, so does the whole theory. 

The language of our theory will be 

$$F \cup \{ +, -, \cdot, ^{-1}, O, < \} \cup \{  c_{z_i} \mid i \ge - 1   \}  $$
where $F$ is the set of all constants $c_{a}$ for each $a \in D$, $+$, $\cdot$ and $<$ are binary function symbols, 
$-$ and $^{-1}$ are unary function symbols, $O$ is an unary relation symbol and $c_{z_i}$ is a constant for all $i \ge 0$. 
Let $\mathcal{A}$ be the theory of the quotient division ring of the ring  $R[y ; \delta]$. By $\mathcal{B}$ we will denote the set of axioms for valuation rings $O$ on division rings:
\begin{enumerate}
\item $O(0) \wedge O(1)$
\item $\forall a : O(a) \vee O(a^{-1})$
\item  $\forall a, b  : O(a) \wedge O(b) \rightarrow O(a + b) \wedge O(ab) \wedge O(ba) $
\end{enumerate}

We add all sentences $\mathcal{C}$ that will give proper meaning to the constants $c_{z_i}$ for all $i \ge -1$:
\begin{enumerate}
\setcounter{enumi}{3}
\item $c_{ z_{0}} = c_{y}  $
\item $c_{z_{i + 1}} = c_{z_{i  }} - c_{\gamma_{i + 1} x^{-r_{i + 1}}   }         $   
\item $O(c_{x^{r_{i + 1}}z_i  }) \wedge O(c_{(x^{r_{i + 1}}z_i)^{-1}  })  $
\end{enumerate}

Our theory is then the union of all the above axioms from $\mathcal{A}$ to $\mathcal{C}$. Since all finite subsets of the theory have a model, namely, the valuation $v_n$ described in the beginning of this subsection, so does, by compactness, the whole theory. Since the theory contains $F$, the set of all constants $c_a$ for each $a \in D$, the models are valued division rings which all contain $D$. We pick a model of the theory, a pair $(D_1, v)$, where $D_1$ is a division ring with valuation $v$.

We now show that the $v$-value  is uniquely determined for every $f \in R[y ; \delta]$. It will then follow that $v$ is uniquely determined on the whole quotient ring $D$.
Every $f \in R[ y ; \delta]$ can be written as
$$
f = \sum_{i = 0}^n p_i^{(0)}(x) y^i, \text{  with $p_i^{(0) }(x) \in R$ for each $0 \le i \le n$}.
$$

For the time being, we ignore the terms we get when we change the order of multiplication. At the end of this subsection, we will see that they do not influence $v(f)$. For each $k \ge 1$, we define
$$ a_k := y - z_k = \sum_{i = 1}^k \gamma_i x^{- r_i}
$$
and write
\begin{align*}
f = \sum_{i = 0}^n p_i^{(0)}(x) y^i = \sum_{i = 0}^n p_i^{(0)}(x) (a_k + z_k  ) ^i \sim \sum_{i = 0}^n p_i^{(0)}(x) \sum_{j = 0}^i  {i \choose j }a_{k}^{i - j} z_k^j = \sum_{j = 0}^n p_j^{(k)}(x)z_k^j   
\end{align*}
with $$p_j^{(k)}(x) := \sum_{i = j}^n  {i \choose j } p_i^{(0)}(x)   a_{k}^{i - j}  $$ for each $0 \le j \le n$ and $k \ge 1$. For each $0 \le j \le n$,
$g_j(t) := \sum_{i = j}^n {i  \choose j  } p_i^{(0)}(x)t^i  $ is a polynomial in $R[t]$. The quotient field of
$$C := \{ \sum_{k \ge m} a_k x^{- \frac{k}{n} }  \mid a_k \in \C, m \in \Z, n \in \N  \}   $$
is the algebraic closure of the quotient field of $R$, as shown in for example \cite{MARK}. Since $\sum_{i = 1}^{\infty}\gamma_i x^{- r_i} $ is not in the quotient field of $C$, it is not a root of $g_j(t)$ for any $0 \le j \le n$. We conclude that for some $K \in \N$,
$v( p_j^{(k)}(x)   )   = v ( p_j^{(K)}(x)  ) $ for all $ k \ge K  $  and all $0 \le j \le n$.
We can then write 
\begin{align*}
f = \sum_{i = 0}^n p_i^{(0)}(x) y^i = \sum_{i = 0}^n p_i^{(K)}(x) z_K^i =  \sum_{i = 0}^n   p_i^{(k)}(x) z_k^i
\end{align*}
with $v(p_i^{(k)}(x)  ) = v(p_i^{(K)}(x))$ for all $k \ge K$. 

We now show that from some $K' \ge K$,  $v(p_0^{(k)}   ) < v(p_i^{(k)} z_{k}^i    )  $ for all $1 \le i \le n$ and all $k > K'$. For all $k > K$,
$$p_0^{(k + 1)}(x) = \sum_{i = 0}^n p_i^{(k)}(x)(a_{k + 1} - a_{k})  ^i =  \sum_{i = 0}^n p_i^{(k )}(x)(\gamma_{k + 1} x^{ - r_{k + 1}  })  ^i  $$
with $v(  p_i^{(k )}(x)   ) = v(  p_i^{(K)}  (x) )  $. Since $(r_i)_{i \ge 1}$ is an increasing sequence with  $\lim_{i  \rightarrow \infty} r_i = 1$, there exists $ K' \ge K $ such that 
$$ \min_{i = 0, \ldots, n} \{ v( p_i^{(K')}(x) )  + i r_{K' + 1}    \} = \min_{i = 0, \ldots, n} \{ v( p_i^{(K)}(x) )  + i r_{K'+ 1}    \}  $$ is achieved at exactly one $0 \le i \le n$. We conclude $v( p_0^{(K' +  1)}   ) = v( p_0^{(K)}   ) $. Then for each $1 \le i \le n$,
$$ v( p_i^{(k)}(x) z_k^i  ) = v( p_i^{(k)}(x)  )  + i r_{k+ 1}  =   v( p_i^{(K)}(x)  )  + i r_{k+ 1} >   v( p_i^{(K)}(x)  )  + i r_{k} \ge v( p_0^{(K)}   ).$$

To show that $v(f)$ is equal to $v(p_0^{(K)}(x)) \in \Q$, we must show that the $v$-value of all the terms we get when we change the order of multiplication must be strictly greater than $v(p_0^{(K)}(x))$. For all $k \ge 1$, we write
\begin{align*}
p_0^{(k)}(x) = \sum_{i = 0}^n p_i^{(0)}(x) a_k^i = \sum_{i = 0}^n p_i^{(0)}(x) ( \sum_{j = 1}^k \gamma_j x^{ - r_j } )^i
= \sum_{i = 0}^n  \sum_{j = 1}^{k_i}  p_i^{(0)}(x) \alpha_i x^{- q_j},
\end{align*}
with $q_j \in \Q$ for all $1 \le j \le k_i$ and $1 \le i \le n$.
Since for all $0  \le i \le n$,  $p_i^{(0)}(x) \in R$ and $(r_i)_{i \ge 1}$ is an increasing sequence with  $\lim_{i  \rightarrow \infty} r_i = 1$,
there exists some $k \ge 1$ such that the term of the sum with $v$-value
$$ \min_{i = 1, \ldots, n  } \{ v( p_i^{(0)}(x)) + (i - 1) r_1 + r_k    \}  $$
is the only term in the sum with its $v$-value. We conclude
\begin{align}\label{basevalue}
v(p_0^{(K)}(x) ) = v(p_0^{(k)}(x) )   \le v( p_i^{(0)}(x)) + (i - 1) r_1 + r_k 
\end{align}
for all $1 \le i \le n$. On the other hand, we can write
\begin{align*}
f = \sum_{i = 0}^n p_i^{(0)}(x) y^i = \sum_{i = 0}^n p_i^{(0)}(x)  (a_k + z_k)^i.
\end{align*}
For all $0 \le i \le n$, all terms of $(a_k + z_k)^i$ when expaned are of the form $a_k^{\ell_1} z_k^{ \ell_2} \ldots a_k^{\ell_{i - 1} }z_k^{ \ell_{i}}$ with $\ell_1, \ldots, \ell_i \ge 0$ and $\ell_1 + \cdots + \ell_i = i$. Since $v(a_k ) = r_1$, $v(z_k) = r_{k + 1}$ and 
$$v[a_k, z_k] = v[ \sum_{j = 1}^k \gamma_j x^{ - r_j }, y - \sum_{j = 1}^k \gamma_j x^{ - r_j }  ]  = v( \sum_{j = 1}^k \gamma_j [x^{ - r_j }, y]  ) = 1 + r_1,$$
it follows that the $v$-value of each term we get when we change the order of multiplication is at least $v(  p_i^{(0)}(x)  ) + (i - 1)r_1 + 1$ for each $1 \le i \le n$, which is, as is immediate from (\ref{basevalue}), strictly greater than $v(p_0^{(K)}(x) )$. We thus conclude
$v(f) = v(p_0^{(K)}(x)) \in \Q$. It also follows that $v(f) = v_k(f)$ for all $k \ge K'$ with $v_k$ as defined in the beginning of this section.
As every element of $D$, the quotient ring for $R[y ; \delta]$, can be written as $fg^{-1}$ with $f, g \in R[y ; \delta]$, it follows that $v$ is uniquely determined on $D$. 
We see that the value group for $v$ is equal to $\Q$. We conclude from Corollary \ref{ratvalgrp} that $v$ is strongly abelian.

It remains to show that the residue field for $v$ is equal to $\R$. Suppose $v(f) = 0$ for some $f \in D$. Then $v_k(f) = 0$ for all $k \ge K$ for some $K \in \N$. We can write
\begin{align*}
f = { (\sum_{i = 0}^m p_i^{(K)} z_K^i ) }{ ( \sum_{j = 0}^n q_j^{(K)} z_K^j )^{-1} }
\end{align*} 
with $ p_i^{(K)},  q_j^{(K)} \in R$ and 
 $$v(   \sum_{i = 0}^m p_i^{(K)} z_K^i   ) = v(  p_0^{(K)}  ) = v(  \sum_{j = 0}^n q_j^{(K)} z_K^j) = v( q_0^{(K)}) = q \in \Q. $$ 
It follows that $v( x^q \sum_{i = 0}^m  p_i^{(K)} z_K^i   ) = v(  x^q\sum_{j = 0}^n q_j^{(K)} z_K^j) = 0$ and
$\alpha := \overline{ x^q \sum_{i = 0}^m  p_i^{(K)} z_K^i    }  = \overline{ x^q p_0^{(K)} } \in \overline{R} = \R$,
$\beta := \overline{ x^q \sum_{j = 0}^n  q_j^{(K)} z_K^j    }  = \overline{ x^q q_0^{(K)} } \in \overline{R} = \R$.
We conclude $\overline{f} = \alpha\beta^{-1} \in \R$. So the residue field for $v$ is indeed equal to $\R$.

\subsection{Extensions of valuations from $\rweyl$ to $R[x ; \delta]$ }

In this section, we characterize the valuations on $\rweyl$ with residue field $\R$ that have an extension to 
$R[y ; \delta]$ with the same residue field. 

Since $x^{ \frac{m_i}{n_i} } \in R[y ; \delta]$, it follows that any valuation $v'$ that extends $v$ to a valuation on $R[y ; \delta]$ with residue field $\R$ must satisfy  $v'(x^{ \frac{m_i}{n_i} }\omega_{i - 1}) = 0$
and $ \tilde{ \gamma_i} := \overline{ x^{ \frac{m_i}{n_i} }\omega_{i - 1}  }  \in \R $. 
In the next proposition, we show the necessary condition for a valuation $v$ on $\rweyl$ to have an extension $v'$ to $\R[y ; \delta]$ with the same residue field. 

\begin{trditev}\label{extendreal}
Let $v$ be a valuation on $\rweyl$ with residue field $\R$ associated to a sequence $( \omega_i )_{i \ge -1}$ with $v( \omega_{i-1}) = \frac{m_i}{n_i}$ for $i \ge 1$ and $\overline{x^{m_i} \omega_{i-1}^{n_i} } = \beta_i \in \R$. Let $\alpha_{i, j} \in \R$ be as in Section \ref{rweyl}. Let $2^{h_i}$ be the greatest power of two dividing $n_i$ for all $i \ge 0$. Then $v$ can be extended to a valuation on $R[y ; \delta]$ with the same residue field only if it fulfils the following conditions:
\begin{enumerate}
\item For each $i$ such that $n_i$ is even, $\beta_i > 0$ must hold, and
\item for each $i, j , \ell $ with $h_i < h_j \le h_{\ell}$, $\alpha_{i, j } \alpha_{i, \ell} > 0$ must hold.
\end{enumerate}
\end{trditev}
\proof
Since 
$$  \tilde{\gamma_i}^{n_i } = \overline{ (x^{ \frac{m_i}{n_i} }\omega_{i - 1} )^{n_i}} = \overline{ x^{m_i}\omega_{i-1}^{n_i}  }   = \beta_i$$
due to Proposition \ref{ratdepelt},  it is obvious that $\tilde{\gamma_i}$ must be equal to an $n_i$-th root of $\beta_i$. If $n_i$ is odd, $\tilde{\gamma_i} \in \R$ is uniquely determined regardless of $\sgn( \beta_i )$, while if $n_i$ is even, $\tilde{\gamma_i} \in \R$ only if $\beta_i > 0$. It is thus obvious that $\beta_i > 0$ must hold for all $i$ where $n_i$ is even if $v$ can be extended from a valuation on $\rweyl$ to a valuation on $R[y ; \delta]$ with the same residue field. This proves the necessity of the first condition.

 To prove the necessity of the second condition, we first observe that
\begin{align*}
\alpha_{i, j } = \overline{ \omega_{i-1}^{K_{i,j}} \omega_{j - 1}^{- K_{j, i}}  } = \tilde{\gamma_i}^{K_{i,j}} \tilde{ \gamma_j }^{- K_{j,i}},
\end{align*}
holds for all $i, j \ge 0$. If $h_i < h_j$, 
$K_{i, j}$ is odd while $K_{j, i}$ is even, so $\sgn( \gamma_i ) = \sgn(\alpha_{i, j }) $ must hold.
We can therefore see that
\begin{align*}
\alpha_{i, j } \alpha_{i, \ell} = 
{\tilde{\gamma_i}^{K_{i,j} + K_{i,\ell}  } \tilde{\gamma_j}^{- K_{j,i}  } \tilde{\gamma_{\ell}}^{- K_{\ell, i}  }}
\end{align*}
for all $i, j, \ell \ge 0$. If $h_i  < h_j \le h_{\ell}$, 
both $K_{i,j}$ and $K_{i,\ell}$ are odd while $K_{j,i}$ and $K_{\ell, i}$ are even. It follows that if $v$ can be extended to a valuation on 
$R[y ; \delta]$ with residue field $\R$, $\alpha_{i, j } \alpha_{i, \ell} > 0$ must hold for all $i, j, \ell \ge 0$ with $h_i  < h_j \le h_{\ell}$.
\endproof

In this section, we show that the conditions $(1)$ and $(2)$ of Proposition \ref{extendreal} are also sufficient for $v$ to have an extension to $R[y ; \delta]$ with residue field $\R$.
Let $v$ be any valuation on $\rweyl$ satisfying the conditions described in Proposition \ref{extendreal}. We will first determine $\tilde{ \gamma_i } \in \R$ for all $i \ge 0$. 
If $n_i$ is odd, there is a unique choice of $\tilde{ \gamma_i } \in \R$. Suppose then $n_i$ is even and $\beta_i > 0$ for some $i \ge 0$. If $h_i < h_j$ for some $j$, then $\sgn(\tilde{\gamma_i}) = \sgn( \alpha_{i,j}  ) = \sgn(\tilde{\gamma_i}^{K_{i,j}} \tilde{ \gamma_j }^{- K_{j,i}})  $ since $K_{i,j}$ is even 
while $K_{j, i}$ is odd. 

We can conclude that if for every power of two $2^{h}$ there is an $i \ge 0$ (or, equivalently, if the $v$-value group is $2$-divisible), $\tilde{\gamma_i} \in \R$ is uniquely determined for all $i \ge 0$. If on the other hand, the value group is non-$2$-divisible, there is an $i \ge 0$ such that $2^{h_i}$ is maximal for all $i \ge 0$. We then have two choices for $\tilde{\gamma_i}$ - a positive or a negative one. The sign of $\tilde{\gamma_j}$ is then uniquely determined for all $j \ge 0$ since
\begin{align*}
\alpha_{i, j } = \tilde{\gamma_i}^{K_{i,j}} \tilde{ \gamma_j }^{- K_{j,i}},
\end{align*}
where $K_{i,j}$ is even and $K_{j, i}$ is odd.
We will now take an arbitrary valuation $v$ on $\rweyl$ satisfying the conditions of Proposition \ref{extendreal}, constructed by a sequence of $( \omega_i )_{i \ge -1}$ as shown in Section \ref{rweyl}. We also pick $\tilde{ \gamma_i} \in \R$ for all $i$.
Then we will describe $v$'s extension to $R[y ; \delta]$ with a sequence of  $( z_i )_{i \ge 0}$ like in the beginning of Section \ref{middle}.

Suppose the valuation $v$ on $\rweyl$ is given by a sequence  $(\omega_i)_{i \ge -1}$  with $\omega_{-1} = x, \omega_0 = y$ and $\omega_i = x^{m_i} \omega_{i-1}^{n_i} - \beta_i$, $\beta_i \in \R$ and $\gcd( m_i, n_i ) =1$ for all $i$. Suppose also that $v$ satisfies conditions $(1)$ and $(2)$ of Proposition \ref{extendreal}. We will show that there is exactly one extension of $v$ from $\rweyl$ to $R[y ; \delta]$ for each appropriate choice of 
$( \tilde{\gamma_i}  )_{i \ge 1}$.

\begin{lema} \label{written form}
For each $k \ge \ell \ge 0$, $\omega_k$ can be written in the following form:
\begin{align*}
\omega_k &= ( \Pi_{i = \ell + 1}^k x^{ \frac{m_i}{n_i}  })  \omega_{\ell} (  \Pi_{i =  \ell + 1}^k B_i    )   - \sum_{i = \ell + 1}^k ( \Pi_{j = i + 1}^k x^{ \frac{m_j}{n_j}  }) \tilde{\gamma_i} ( \Pi_{j = i}^k B_j  ) \\
 &+ \sum_{i = \ell + 1}^k ( \Pi_{j = i + 1}^k x^{ \frac{m_j}{n_j} } ) A_i  ( \Pi_{j = i + 1}^k B_i      ),
\end{align*}
with
\begin{align*}
A_i &= \sum_{j = 1}^{n_i - 1} ( x^{ \frac{m_i}{n_i} }  \omega_{i - 1}  )^{n_i - j - 1} x^{ \frac{m_i}{n_i}  } [ x^{ j\frac{m_i}{n_i}  }  , \omega_{i - 1}   ] \omega_{i - 1}^j \\
B_i &= \sum_{j = 1}^{n_j} ( x^{ \frac{m_i}{n_i} }  \omega_{i - 1}  )^{n_i - j} \tilde{\gamma_i}^{j - 1} \\
\end{align*}
for each $1 \le i \le k$. 
\end{lema}
\proof
We prove the lemma by induction on $k \ge \ell$. For $k = \ell$, it is trivially true since we get $\omega_{k} = \omega_k$.

Suppose now the equation holds for some $k \ge \ell$. Then 
\begin{align*}
\omega_{k + 1} &= x^{m_{k + 1}} \omega_k^{n_{k + 1}} - \beta_{k + 1} = 
( x^{ \frac{m_{k + 1}}{ n_{k + 1  }  } } \omega_k - \tilde{\gamma}_{k + 1}    )B_{k + 1} + A_{k + 1} \\
&=  
x^{ \frac{m_{k+1}}{ n_{k+1} } }( ( \Pi_{i = \ell + 1}^k x^{ \frac{m_i}{n_i}  })  \omega_{\ell} (  \Pi_{i =  \ell + 1}^k B_i    )   - \sum_{i = \ell + 1}^k ( \Pi_{j = i + 1}^k x^{ \frac{m_j}{n_j}  }) \tilde{\gamma_i} ( \Pi_{j = i}^k B_j  ) \\
 &+ \sum_{i = \ell + 1}^k ( \Pi_{j = i + 1}^k x^{ \frac{m_j}{n_j} } ) A_i  ( \Pi_{j = i + 1}^k B_i      )    
       )B_{k + 1}
 -  \tilde{\gamma}_{k + 1} B_{k + 1} + A_{k + 1} \\
&= ( \Pi_{i = \ell + 1}^{k + 1} x^{ \frac{m_i}{n_i}  })  \omega_{\ell} (  \Pi_{i =  \ell + 1}^{k + 1} B_i    )   - \sum_{i = \ell + 1}^{k + 1} ( \Pi_{j = i + 1}^{k + 1} x^{ \frac{m_j}{n_j}  }) \tilde{\gamma_i} ( \Pi_{j = i}^{k + 1} B_j  ) \\
 &+ \sum_{i = \ell + 1}^{k + 1} ( \Pi_{j = i + 1}^{k + 1} x^{ \frac{m_j}{n_j} } ) A_i  ( \Pi_{j = i + 1}^{k+ 1} B_i      ) ,
\end{align*}
where we used the induction hypothesis, which is
\begin{align*}
\omega_k &= ( \Pi_{i = \ell + 1}^k x^{ \frac{m_i}{n_i}  })  \omega_{\ell} (  \Pi_{i =  \ell + 1}^k B_i    )   - \sum_{i = \ell + 1}^k ( \Pi_{j = i + 1}^k x^{ \frac{m_j}{n_j}  }) \tilde{\gamma_i} ( \Pi_{j = i}^k B_j  ) \\
 &+ \sum_{i = \ell + 1}^k ( \Pi_{j = i + 1}^k x^{ \frac{m_j}{n_j} } ) A_i  ( \Pi_{j = i + 1}^k B_i      )
\end{align*}
 in the second equation.
\endproof
 Since $v(\omega_{i-1}) = \frac{m_i}{n_i}$ and $\overline{ x^{ \frac{m_i}{n_i} }  \omega_{i - 1}} = \tilde{\gamma_i}$, $v(B_i) = 0$ and
$\overline{B_i} = n_i \tilde{\gamma_i}^{n_i - 1}$ must hold for all $i \ge 1$. From 
\begin{align*}
[x^{m_i}, \omega_{i - 1}] = \sum_{j = 0}^{n_i - 1} x^{ \frac{m_i}{n_i} \cdot j} [x^{\frac{m_i}{n_i}}, \omega_{i - 1} ] x^{ \frac{m_i}{n_i} \cdot ( n_i - j ) },
\end{align*}
we can see that $v(A_i) = v[x^{\frac{m_i}{n_i}}, \omega_{i - 1}] = v[x, \omega_{i-1}] + 1 - \frac{m_i}{n_i} = -v (xy \omega_1 \ldots \omega_{i-1}).  $ 
Since $\sum_{i \ge -1}^k v(\omega_i) < 0$ for every $k$ by Lemma \ref{kompogoj}, $$v ( \sum_{i = \ell + 1}^k( \Pi_{j = i + 1} x^{ \frac{m_j}{n_j} } ) A_i  ( \Pi_{j = i}^k B_i      ) ) > v(\omega_k)$$ will hold for every $k$, which is why we can ignore the terms containing $A_i$ during our evaluations  of $v(z_i)$ where for each $i$, $z_i \in R[y ; \delta]$ is as in the beginning of this section. 


\endproof 
\begin{lema} \label{construction one}
Suppose $v$ is a valuation on $\rweyl$ constructed from a sequence $( \omega_i  )_{i  \ge -1}$ that extends to a valuation on $R[y ; \delta]$ and thus $D$, its quotient division ring. For a given $i \ge 1$, define a sequence $( S_{i, j}   )_{j \ge 1}$ by:
\begin{enumerate}
\item[(a)] $S_{i, 1} := (x^{ \frac{m_i}{n_i}}\omega_{i-1} - \tilde{\gamma_i})^{-1 }(B_i - \overline{B_i}  )$,
\item[(b)]$S_{i, j+1} := (x^{ \frac{m_i}{n_i}}\omega_{i-1} - \tilde{\gamma_i})^{-1 }(S_{i,j} - \overline{ S_{i, j}  }   )$ for $j \ge 1$.
\end{enumerate}
Then for each $j \ge 1$:
\begin{enumerate}
\item $S_{i, j } = \sum_{k = 1}^{n_i - j} N_{k, j} ( x^{ \frac{m_i}{n_i}}\omega_{i-1})^{n_i - j - k} \tilde{\gamma_i}^{k - 1}   $ with $N_{k, 1} = k$ and 
$N_{k, j+1} = \sum_{\ell = k}^{n_i - j}N_{\ell,j}$, and
\item $v( S_{i, j}  ) = 0$, $\overline{S_{i, j}} = N_{1, j +1} \tilde{\gamma_i}^{n_i - j - 1}  $
\end{enumerate}
for all $1 \le j \le n_i- 1$. 
\end{lema}

\proof
We prove the first statement of the lemma by induction on $j \ge 1$. To show the basis of induction, we evaluate
\begin{align*}
B_i - \overline{B_i} &= (x^{ \frac{m_i}{n_i}  }\omega_{i-1}  - \tilde{\gamma_i}   )( \sum_{k = 1}^{n_i} ( (x^{ \frac{m_i}{n_i}    } \omega_{i - 1}   )^{n_i - k - 1} + (x^{ \frac{m_i}{n_i}    } \omega_{i - 1}   )^{n_i - k - 2}\tilde{\gamma_i}  +  \cdots  
+  \tilde{\gamma_i}^{n_i - k- 1}   ) \tilde{\gamma_i}^{k - 1}    ) \\
 &=  (x^{ \frac{m_i}{n_i}  }\omega_{i-1}  - \tilde{\gamma_i}   ) \sum_{k = 1}^{n_i - 1} k  (x^{ \frac{m_i}{n_i}    } \omega_{i - 1}    )^{n_i  - k - 1} \tilde{\gamma_i}^{ k - 1}.
\end{align*}
We can thus see  $S_{i, 1} = \sum_{k = 1}^{n_i - 1} k  (x^{ \frac{m_i}{n_i}    } \omega_{i - 1}    )^{n_i  - k - 1} \tilde{\gamma_i}^{ k - 1}$ and since $\overline{x^{ \frac{m_i}{n_i}    } \omega_{i - 1}} = \tilde{\gamma_i}$,
we see that the lemma holds in case $j = 1$.
Now we suppose that the statement is true for some $j \ge 1$, i.e., $S_{i, j } = \sum_{k = 1}^{n_i - j} N_{k, j} ( x^{ \frac{m_i}{n_i}}\omega_{i-1})^{n_i - j - k} \tilde{\gamma_i}^{k - 1}   $ and $v( S_{i, j}  ) = 0$, $\overline{S_{i, j}} = N_{1, j +1} \tilde{\gamma_i}^{n_i - j - 1}  $. We then write
\begin{align*}
S_{i, j } - \overline{S_{i, j}} &= \sum_{k = 1}^{n_i - j} N_{k, j} ( x^{ \frac{m_i}{n_i}}\omega_{i-1})^{n_i - j - k} \tilde{\gamma_i}^{k - 1} - N_{1, j +1} \tilde{\gamma_i}^{n_i - j - 1} \\
&=  \sum_{k = 1}^{n_i - j}N_{k, j} ( x^{ \frac{m_i}{n_i}}\omega_{i-1})^{n_i - j - k} \tilde{\gamma_i}^{k - 1} -   \sum_{k = 1}^{n_i - j} N_{k, j}\tilde{\gamma_i}^{n_i - j - 1} \\
&= \sum_{k = 1}^{n_i - j - 1}  N_{k, j}(  ( x^{ \frac{m_i}{n_i}}\omega_{i-1})^{n_i - j - k} -  \tilde{\gamma_i}^{n_i - j - 1}     ) \tilde{\gamma_i}^{k - 1} \\
&= (x^{ \frac{m_i}{n_i}  }\omega_{i-1}  - \tilde{\gamma_i}   )\sum_{k = 1}^{n_i - j  - 1}  N_{k, j} (  ( x^{ \frac{m_i}{n_i}}\omega_{i-1})^{n_i - j - 1 - k}     + \cdots       
+ \tilde{\gamma_i}^{n_i - j - 1 - k}
)\tilde{\gamma_i}^{k - 1} \\
&=  (x^{ \frac{m_i}{n_i}  }\omega_{i-1}  - \tilde{\gamma_i}   )  \sum_{k = 1}^{n_i - j - 1}    (  N_{k, j} + N_{k + 1, j} + \cdots 
+ N_{ n_i - j - 1, j   }) ( x^{ \frac{m_i}{n_i}}\omega_{i-1})^{n_i - j - k} \tilde{\gamma_i}^{k - 1} \\
&=  (x^{ \frac{m_i}{n_i}  }\omega_{i-1}  - \tilde{\gamma_i}   )  \sum_{k = 1}^{n_i - j - 1}    N_{k, j+ 1} ( x^{ \frac{m_i}{n_i}}\omega_{i-1})^{n_i - j - k} \tilde{\gamma_i}^{k - 1}
\end{align*}
proving the first statement of the lemma. The second statement immediately follows from the first since $\overline{x^{ \frac{m_i}{n_i}  }\omega_{i-1}} = \tilde{\gamma_i}$  and thus 
\begin{align*}
v( S_{i, j}   ) &= v(   \sum_{k = 1}^{n_i - j }   N_{k, j} ( x^{ \frac{m_i}{n_i}}\omega_{i-1})^{n_i - j - k} \tilde{\gamma_i}^{k - 1}  ) = 0, \\
\overline{S_{i, j}  } &=  \overline{  \sum_{k = 1}^{n_i - j}    N_{k, j} ( x^{ \frac{m_i}{n_i}}\omega_{i-1})^{n_i - j - k} \tilde{\gamma_i}^{k - 1}    }  = N_{1, j + 1} \tilde{\gamma_i}^{n_i - j  - 1}.
\end{align*}
\endproof

\begin{lema} \label{construction two}
Supose $v$ is as in Lemma \ref{construction one}. For a given $i \ge 1$, define a sequence $( D_{i, j}   )_{j \ge 1}$ by:
\begin{enumerate}
\item[(a)] $D_{i, 1} = B_1B_2 \cdots B_i - \overline{  B_1B_2 \cdots B_i  }$,
\item[(b)] $D_{i, j+1} = x^{v(D_{i, j})}D_{i, j} - \overline{x^{v(D_{i, j})}D_{i, j}}$.
\end{enumerate}
Then for each $j \ge 1$:
\begin{enumerate}
\item 
 $D_{i, j}$ is a $\R$-linear sum of terms which are products of elements from the set
\begin{align}\label{factors}
 \{  \omega_i     \}_i \cup \{  B_i  \}_i \cup \{ B_i^{-1}  \}_i \cup  \{ S_{i,j}  \}_{i, j }, 
\end{align}
where parts of the product are conjugated by a rational power of $x$.
\item $v( D_{i,j}  )$ is a sum of $v( \omega_{\ell})$ for finitely many $\omega_{ \ell }$.
\item $\overline{x^{v(D_{i, j})}D_{i, j}}$ is sum of products of $\tilde{\gamma_k}$ for various $k$.
\end{enumerate}
\end{lema}
\proof
Since
\begin{align*}
B_1B_2 \cdots B_i  - \overline{ B_1B_2 \cdots B_i } &= \sum_{j = 1}^i B_1 \cdots B_{j - 1}( B_j - \overline{B_j}   ) \overline{ B_{j + 1} \cdots B_i} \\
&=  \sum_{j= 1}^i  B_1 \cdots B_{j - 1}\omega_j B_{j}^{-1} S_{j, 1} \overline{ B_{i + 1} \cdots B_i},
\end{align*}
and 
\begin{align*}
x^{ \frac{m_{j + 1} }{n_{j + 1}}  }  B_1 \cdots B_{j - 1}\omega_j B_{j}^{-1} S_{j, 1} \overline{ B_{i + 1} \cdots B_i}
=  (x^{ \frac{m_{j + 1} }{n_{j + 1}}  } B_1 \cdots B_{j - 1} x^{ - \frac{m_{j + 1} }{n_{j + 1}}  } )x^{ \frac{m_{j + 1} }{n_{j + 1}}  }
\omega_j B_{j}^{-1} S_{j, 1} \overline{ B_{i + 1} \cdots B_i},
\end{align*} 
 for each $1 \le j \le k$, the first two statements of the lemma follow from:
\begin{enumerate}
\item $x^{\frac{m_i}{n_i}}\omega_{i - 1} - \tilde{\gamma}_i =  \omega_i B_i^{-1}$  ,
\item $B_i - \overline{B_i } = (x^{\frac{m_i}{n_i}}\omega_{i - 1} - \tilde{\gamma})S_{i, 1} =  \omega_i B_i^{-1} S_{i, 1} $,
\item $S_{i, j } - \overline{S_{i,j}} = (x^{\frac{m_i}{n_i}}\omega_{i - 1} - \tilde{\gamma} S_{i, j + 1}  =  \omega_i B_i^{-1}    S_{i, j+1}  $,
\item $B_i^{-1 }- \overline{B_i }^{-1} =  - {B_i}^{-1} (B_i - \overline{B_i })  \overline{B_i}^{-1} =  -
{B_i}^{-1}  \omega_i B_i^{-1} S_{i, 1} \overline{B_i}^{-1}  $,
\end{enumerate}
where we ignore the terms we get when we change the order of multiplication. We can do that that since these terms are procucts of $A_i$ as defined in Lemma \ref{written form} and terms with zero $v$-value. As we have already mentioned, these terms will not influence the construction of the extension of a valuation on $\rweyl$ to $R[y ; \delta]$. 
Indeed - we can see by induction on $j \ge 1$ that each term of the sum is a product of factors equal to, modulo conjugation by a rational power of $x$, one of the elements of  the set (\ref{factors}); that is,
\begin{enumerate}
\item either equal to $\omega_i$, or
\item equal to a power of $B_j  $ or $S_{j, i}$ for some $i,j$. 
\end{enumerate}
Since the latter have $v$-value equal to zero and since both $B_i - \overline{B_i}$ and $S_{i, j} - \overline{ S_{i, j} }$ are products of 
$\omega_i$,  a power of $B_i^{-1}$ and $S_{i, j} $ for some $i, j$, $v(D_{i, j})$ is a sum of $v(\omega_i)$ for some $i$.
The last statement of the lemma follows from the fact that for all $i, j$, $\overline{B_{i}}$ and $\overline{S_{i ,j}}$ are of the form $N \tilde{\gamma}_i$ for $N \in \N$.
\endproof

If $v$ is a valuation on $R[y ; \delta]$ with residue field $\R$, then, as we have presented in Section \ref{middle}, $v$ can be constructed from a sequence $( z_i )_{i \ge 0} \subset R[y ; \delta]$.  
In the next proposition, we make the first comparison between this construction and the construction of a valuation on $\rweyl$ from a sequence $ ( \omega_i)_{i \ge -1}   $ described in Section \ref{rweyl}.

\begin{lema} \label{construction three}
Suppose $v$ is as in Lemma \ref{construction one}. Suppose that for some $k, \ell \ge 0$, we can write
\begin{align*}
\omega_k = (\Pi_{i = 1}^k x^{ \frac{m_i}{n_i} }) z_{\ell} (\Pi_{i = 1}^k B_i )+ C,
\end{align*}
where $C$ is a $\R$-linear sum of elements of the form $D_{i,j}$ for some $i, j\ge 0$. Then:
\begin{enumerate}
\item If $v( \omega_k   ) > v( \Pi_{i = 1}^k x^{ \frac{m_i}{n_i} }    z_{\ell} ) \in \Q,$ then
\begin{align*}
\omega_k = (\Pi_{i = 1}^k x^{ \frac{m_i}{n_i} }) z_{\ell + 1} (\Pi_{i = 1}^k B_i )+ C_1.
\end{align*}
\item 
 If $v( \omega_k   ) =  v( \Pi_{i = 1}^k x^{ \frac{m_i}{n_i} }    z_{\ell} ) \in \Q $, then 
\begin{align*}
\omega_{k + 1}  = (\Pi_{i = 1}^{k + 1} x^{ \frac{m_i}{n_i} }) z_{\ell + 1} (\Pi_{i = 1}^{k + 1} B_i )+ C_2.
\end{align*}
\item 
 If $v( \omega_k   ) <  v( \Pi_{i = 1}^k x^{ \frac{m_i}{n_i} }    z_{\ell} ),$ then, if $v( \omega_k   ) \in \Q$,
\begin{align*}
\omega_{k + 1}  = (\Pi_{i = 1}^{k + 1} x^{ \frac{m_i}{n_i} }) z_{\ell } (\Pi_{i = 1}^{k + 1} B_i )+ C_3.
\end{align*}
\end{enumerate}
Here, $C_1, C_2$ and $C_3$ are other $\R$-linear sums of elements of the form $D_{i, j}$ as in Lemma \ref{construction three} for some $i, j \ge 0$.
\end{lema}

\proof
Suppose first $v( \omega_k   ) > v( \Pi_{i = 1}^k x^{ \frac{m_i}{n_i} }    z_{\ell} ) = v(C)$.
Since $v(C) = r_{\ell + 1} - \sum_{i = 1}^k \frac{m_i}{n_i} \in \Q$, then $r_{\ell + 1} := v(z_{\ell}) \in \Q$ as well.
So, for $z_{\ell + 1} = z_{\ell} - x^{- r_{\ell + 1}}\gamma_{\ell + 1}$, we can write
\begin{align*}
\omega_k &= (\Pi_{i = 1}^k x^{ \frac{m_i}{n_i} }) z_{\ell + 1} (\Pi_{i = 1}^k B_i ) +  \gamma_{\ell + 1}(\Pi_{i = 1}^k x^{ \frac{m_i}{n_i} })x^{- r_{\ell + 1}} (\Pi_{i = 1}^k B_i )    + C \\
&=  (\Pi_{i = 1}^k x^{ \frac{m_i}{n_i} }) z_{\ell + 1} (\Pi_{i = 1}^k B_i ) +  x^{- v(C)}( \gamma_{\ell + 1} (\Pi_{i = 1}^k B_i )     + x^{v(C)}C)
\end{align*}
Since $v( z_{\ell + 1}  ) > v( z_{\ell}  )$, we see that $v( \gamma_{\ell + 1} \Pi_{i = 1}^k B_i  +  x^{ v(C)}C)  > 0$. It follows that
$\gamma_{i + 1} = -  \Pi_{i = 1}^k \overline{B_i}^{-1}  \overline{  x^{- v(C)}C}     $, hence
\begin{align*}
C_1 := \gamma_{\ell + 1} \Pi_{i = 1}^k B_i  +  x^{- v(C)}C &= \gamma_{\ell + 1} \Pi_{i = 1}^k B_i +  \overline{  x^{- v(C)}C} - \overline{  x^{- v(C)}C}  +    x^{- v(C)}C \\
&=   \gamma_{\ell + 1}( \Pi_{i = 1}^k B_i  -  \Pi_{i = 1}^k \overline{B_i}  ) + ( x^{- v(C)}C  - \overline{x^{- v(C)}C }   )
\end{align*}
We see that since $C$ is an $\R$-linear sum of $D_{i ,j}$ for various $i,j$, so is, by Lemma \ref{construction two}, $ x^{- v(C)}C  - \overline{x^{- v(C)}C }$. Hence, $C_1$ is an $\R$-linear sum of $D_{i, j}$. 

Now consider the case $v( \omega_k   ) = v( \Pi_{i = 1}^k x^{ \frac{m_i}{n_i} }    z_{\ell} ) \le v(C)$. Since 
$v( \omega_k  ) =  \frac{m_{k + 1}}{ n_{k + 1}  } \in \Q,$ we can evaluate
\begin{align*} 
\tilde{\gamma_{k + 1}} &= \overline{ x^{  \frac{m_{k + 1}}{ n_{k + 1}  }  } \omega_k    } = \overline{ (\Pi_{i = 1}^{k + 1} x^{ \frac{m_i}{n_i} }) z_{\ell} (\Pi_{i = 1}^{k  } B_i ) }   + \overline{ x^{\frac{m_{k + 1} }{n_{k + 1}}   } C } \\
&= \gamma_{\ell + 1} \Pi_{i = 1}^{k  } \overline{B_i }   + \overline{ x^{\frac{m_{k + 1} }{n_{k + 1}}   } C },
\end{align*}
and then deduce
\begin{align*}
\omega_{k + 1} &=  (x^{  \frac{m_{k + 1}}{ n_{k + 1}  }  } \omega_k -  \tilde{\gamma_{k + 1}})B_{k + 1} \\
&=   (\Pi_{i = 1}^{k + 1} x^{ \frac{m_i}{n_i} }) z_{\ell} (\Pi_{i = 1}^{k +1 } B_i ) +  x^{\frac{m_{k + 1} }{n_{k + 1}}   } CB_{k + 1} - \tilde{\gamma_{k + 1}}B_{k + 1} \\
&= (\Pi_{i = 1}^{k + 1} x^{ \frac{m_i}{n_i} }) z_{\ell + 1}  (\Pi_{i = 1}^{k +1 } B_i ) + \gamma_{ \ell + 1}(\Pi_{i = 1}^{k +1 } B_i )  +  x^{\frac{m_{k + 1} }{n_{k + 1}}   } CB_{k + 1} - \tilde{\gamma_{k + 1}}B_{k + 1} \\
&=  (\Pi_{i = 1}^{k + 1} x^{ \frac{m_i}{n_i} }) z_{\ell + 1}  (\Pi_{i = 1}^{k +1 } B_i ) + \gamma_{\ell + 1}( \Pi_{i = 1}^{k  } B_i   - 
 \Pi_{i = 1}^{k  } \overline{B_i}  )B_{k + 1} + ( { x^{\frac{m_{k + 1} }{n_{k + 1}}   } C } - \overline{ x^{\frac{m_{k + 1} }{n_{k + 1}}   } C }  )B_{k + 1} \\
&= (\Pi_{i = 1}^{k + 1} x^{ \frac{m_i}{n_i} }) z_{\ell + 1}  (\Pi_{i = 1}^{k +1 } B_i ) + C_2 
\end{align*}
where $C_2$ is, again, an $\R$-linear sum of $D_{i, j }$ for some $i, j$.

Lastly, we consider the case  $v( \Pi_{i = 1}^k x^{ \frac{m_i}{n_i} }    z_{\ell} ) > v( \omega_k   ) = v(C) \in \Q. $
In this case, $\tilde{\gamma_{k + 1}} =  \overline{ x^{  \frac{m_{k + 1}}{ n_{k + 1}  }  } \omega_k    } =  \overline{ x^{\frac{m_{k + 1} }{n_{k + 1}}   } C } $, so we can write
\begin{align*}
\omega_{k + 1} &=  (x^{  \frac{m_{k + 1}}{ n_{k + 1}  }  } \omega_k -  \tilde{\gamma_{k + 1}})B_{k + 1}  \\
&=  (\Pi_{i = 1}^{k + 1} x^{ \frac{m_i}{n_i} }) z_{\ell }  (\Pi_{i = 1}^{k +1 } B_i ) + ( { x^{\frac{m_{k + 1} }{n_{k + 1}}   } C } - \overline{ x^{\frac{m_{k + 1} }{n_{k + 1}}   } C }  )B_{k + 1}
\end{align*}
and the statement again follows.
\endproof

\begin{izrek} \label{construction result}
Suppose $v$ is a valuation on $\rweyl$, constructed from an either finite or infinite sequence $( \omega_i  )_{i \ge - 1} $ with 
$v( \omega_i) = \frac{m_{i + 1}}{ n_{i + 1} } $. Suppose also that $v$ satisfies the following conditions:
\begin{enumerate}
\item For each $i$ such that $n_i$ is even, $\beta_i > 0$ must hold, and
\item for each $i, j , \ell \ge 0$ with $h_i < h_j \le h_{\ell}$, $\alpha_{i, j } \alpha_{i, \ell} > 0$ must hold.
\end{enumerate}
Then $v$ has a unique extension to a valuation on $R[y ; \delta]$ with residue field $\R$ for each choice of $( \tilde{\gamma_i})_{i \ge 1}  $.
\end{izrek}
\proof
As we know from the beginning of Section \ref{middle}, each valuation on $R[y ; \delta]$ and $D$ with residue field $\R$ can be constructed by either a finite or an infinite sequence $ ( z_i  )_{i \ge 0} $. 
For every sequence $(\omega_i   )_{i \ge - 1}$, we will use the lemmas proved in this section to find the unique sequence  $ ( z_i  )_{i \ge 0} $ which, as we have shown in the beginning of this section, uniquely determines a valuation on $R[y ; \delta]$. Our calculations will then show that the valuation on $D$ defined by the sequence $( z_i )_{i \ge 0}  $ is the extension of the valuation on $\rweyl$ associated to the sequence $(\omega_i)_{i \ge -1}$. 

We determine the finite or infinite sequence $(z_i)_{i \ge 0}$ associated to $v$'s extension to $R[y ; \delta]$.
In the first step, we consider $v(y)$. 
If $v(y) \not\in \Q$, then $v$'s extension to $R[y ; \delta ]$ is clearly uniquely determined, namely the one defined by
\begin{align*}
v(  \sum_{i = 0}^n p_i(x) y^i   ) = \min_{0 \le i \le n}   \{ v( p_i(x)   ) + i v(y)   \}
\end{align*}
for every $ \sum_{i = 0}^n p_i(x) y^i  \in R[ y ; \delta ]$. 

So suppose $v(y) = \frac{m_1}{n_1} \in \Q$. Then, in the second step of our evaluation, we write 
\begin{align*}
\omega_1 = x^{m_1}y^{n_1} - \beta_1 = ( x^{ \frac{m_1}{n_1}}y - \tilde{\gamma_1}   )B_1 = x^{ \frac{m_1}{n_1} } (y -  x^{ - \frac{m_1}{n_1}}\tilde{\gamma_1}   )B_1
\end{align*}
and since $v( \omega_1) > 0$, $v(y -  x^{ - \frac{m_1}{n_1}}\tilde{\gamma_1}) >  \frac{m_1}{n_1} = v(y)$.
We deduce  $z_1 =   y - x^{ - \frac{m_1}{n_1}}{\gamma_1}$ with $\gamma_1 = \tilde{\gamma_1}$.  Obviously, $v(z_1) = v(\omega_1) + \frac{m_1}{n_1} \in \Q$ if and only if $v(\omega_1) \in \Q$. If either and hence both values are irrational, we get a unique extension of $v$ to $R[y ; \delta]$.  
Otherwise, if $v(\omega_1) = \frac{m_2}{n_2} \in \Q$ and hence $r_2 = v(z_2) = \frac{m_1}{n_1} + \frac{m_2}{n_2}$, we continue with the third step of our evaluation by writing
\begin{align*}
\omega_2 = (x^{\frac{m_2}{n_2}} \omega_1  - \tilde{\gamma_2}  ) B_2 =  (x^{\frac{m_2}{n_2}} x^{\frac{m_1}{n_1}}  z_1 B_1  - \tilde{ \gamma_2}) B_2.
\end{align*}
Since $v( \omega_2) > 0$, we conclude that $\gamma_2 = \overline{x^{\frac{m_2}{n_2}} x^{\frac{m_1}{n_1}}  z_1  }$ must be equal to $\tilde{\gamma_2} \overline{B_1}^{-1}$. To evaluate the $v$-value of $z_2 = z_1 - \gamma_2 x^{- r_2}$, we write
\begin{align*}
\omega_2 =  (x^{\frac{m_2}{n_2}} x^{\frac{m_1}{n_1}}  z_2 B_1  +  \gamma_2 B_1   - \tilde{ \gamma_2}) B_2
 = x^{\frac{m_2}{n_2}} x^{\frac{m_1}{n_1}}  z_2 B_1B_2  + (\gamma_2 B_1 - \tilde{\gamma_2}) B_2.
\end{align*}
We note that $\omega_2$ is here written as a sum of $\Pi_{i = 1}^2 x^{  \frac{m_i}{n_i} } z_2 \Pi_{i = 1}^2 B_i + C$ where $C$ is as in Lemma \ref{construction three}.
To determine $v(z_2)$, we compare $v(\omega_2)$ and $v( \gamma_2 B_1 - \tilde{\gamma_2} )$, the latter being equal to $v(\omega_1)$ since
$\gamma_2 B_1 - \tilde{\gamma_2} = \gamma_2(B_1 - \overline{B_1}) $. There are three possible cases:
\begin{enumerate}
\item If $v(\omega_2) < v(\omega_1)$, then $v( x^{\frac{m_2}{n_2}} x^{\frac{m_1}{n_1}} z_2 )$ must be equal to $v( \omega_2)$, so $r_3 = v (z_2)$ is determined by $v(z_2) = v(\omega_2) - \frac{m_1}{n_1} - \frac{m_2}{n_2}$.
If $v(\omega_2) = \frac{m_3}{n_3} \in \Q$, then $v(z_2) = r_3 \in \Q$
and $\gamma_3 = \tilde{\gamma_3}(\overline{B_1B_2})^{-1}$ must hold. It follows that $v( z_2) \in \Q$ if and only if $v(\omega_2) \in \Q$. In this case, the $v$-value of $z_3 = z_2 - \gamma_3 x^{ - r_3  }$ will be determined in the subsequent steps, i.e., by considering $ (v (\omega_i))_{i \ge 3}  $ and
 $ (\tilde{ \gamma_i})_{i \ge 3}    $. By Lemma \ref{construction three}, $\omega_3 = \Pi_{i = 1}^3 x^{ \frac{m_i}{n_i}  } z_3  \Pi_{i = 1}^3 B_i + C_1$, $C_1$ being an $\R$-linear sum of $D_{i, j }$.
\item If $v(\omega_2) = v(\omega_1)$, then $v(z_2)$ depends on 
$\tilde{\gamma_3}$:
\begin{enumerate}
\item If $\tilde{\gamma_3} = \overline{ x^{\frac{m_2}{n_2}} ( \gamma_2 B_1 - \tilde{\gamma_2}  ) B_2}$, then $v( x^{\frac{m_2}{n_2}} x^{\frac{m_1}{n_1}} z_2  )$ must be greater than $v(\omega_1)$ since in this case, $\omega_2 \sim (\gamma_2B_1 - \tilde{\gamma}_2   )B_2  $  and will be determined in the subsequent steps, 
i.e., by considering $ (v (\omega_i))_{i \ge 3}  $ and
 $ (\tilde{ \gamma_i})_{i \ge 3}   $. By Lemma \ref{construction three}, $\omega_3 = \Pi_{i = 1}^3 x^{ \frac{m_i}{n_i}  } z_2  \Pi_{i = 1}^3 B_i + C_2$, $C_2$ being an $\R$-linear sum of $D_{i, j }$.
\item Otherwise,  $v( x^{\frac{m_2}{n_2}} x^{\frac{m_1}{n_1}}  z_2  ) = v(\omega_2)$. In this case, we get $$\gamma_3 = \tilde{\gamma_3} - \overline{ x^{\frac{m_2}{n_2}} ( \gamma_2 B_1 - \tilde{\gamma_2}  )B_2}.$$ By Lemma \ref{construction three}, $\omega_3 = \Pi_{i = 1}^3 x^{ \frac{m_i}{n_i}  } z_3  \Pi_{i = 1}^3 B_i + C_3$, with $C_3$ an $\R$-linear sum of $D_{i, j }$.
\end{enumerate}
\item If $v(\omega_2) > v(\omega_1)$, then $v( x^{\frac{m_2}{n_2}} x^{\frac{m_1}{n_1}}  z_2  ) = v(\omega_1)$ and 
$\gamma_3 = \overline{ - x^{ \frac{m_2}{n_2}  } ( \gamma_2 B_1 - \tilde{\gamma_2}   )B_1^{-1}   }$. By Lemma \ref{construction three}, $\omega_2 = \Pi_{i = 1}^2 x^{ \frac{m_i}{n_i}  } z_3  \Pi_{i = 1}^2 B_i + C_4$, with $C_4$ an $\R$-linear sum of $D_{i, j }$.
\end{enumerate}

The general step of the evaluation is similar to the first three. Suppose that in the previous steps, we have evaluated $v(z_1) = r_2, \ldots, v(z_{\ell - 1}) = r_{\ell}\in \Q$. In the last step, we have, by considering $( \omega_i  )_{i = -1 }^k$ for some $k$, begun to evaluate $v(z_{\ell})$ and we are, with
\begin{align*}
\omega_k &= 
(  \Pi_{i = 1}^k x^{ \frac{m_i}{n_i}  }    ) z_{\ell + 1} ( \Pi_{i = 1}^k B_i   ) + C,
\end{align*}
where $C$ is as in Lemma \ref{construction three}, in one of the five situations:
\begin{enumerate}
\item If $v(\omega_k) < v( C  )$, then $v(z_{\ell}) = v( \omega_k )  - \sum_{i = 1}^k \frac{m_i}{n_i}  $ and 
$\gamma_{\ell + 1} = \tilde{\gamma_{\ell + 1}}\Pi_{i = 1}^k \overline{B_i}^{-1}  $. In case $r_{\ell + 1} = v(z_{\ell} ) \in \Q$, our next step is to evaluate the $v$-value of $z_{\ell + 1} = z_{\ell} - r^{\ell + 1}\gamma_{\ell + 1} $ by writing
\begin{align*}
\omega_{k  + 1 } &= (  \Pi_{i = 1}^{k + 1} x^{ \frac{m_i}{n_i}  }    ) z_{\ell} ( \Pi_{i = 1}^{k + 1} B_i   ) + C_1.
\end{align*}
\item If $v(\omega_k) > v( C  )$, then  $v(z_{\ell}) = v( C  )$ and 
$\gamma_{\ell + 1} = \overline{C}\Pi_{i = 1}^k \overline{B_i}^{-1}  $. In case $r_{\ell + 1} = v(z_{\ell} ) \in \Q$, our next step is to evaluate the $v$-value of $z_{\ell + 1} = z_{\ell} - r^{\ell + 1}\gamma_{\ell + 1} $ by writing
\begin{align*}
\omega_{k   } &= (  \Pi_{i = 1}^{k } x^{ \frac{m_i}{n_i}  }    ) z_{\ell + 1} ( \Pi_{i = 1}^{k } B_i   ) + C_2.
\end{align*}
\item If $v(\omega_k) = v( C  ) \not\in\Q$, then $v(z_{\ell}   ) =  v(\omega_k - C)  - \sum_{i = 1}^k \frac{m_i}{n_i} \not\in \Q$ since in case $ v( \omega_k - C   ) > v (\omega_k)$, $C \sim \omega_k$ must hold, but since, given that $C \neq \omega_k$ and that $C$ is a sum of $D_{i,j}$, $ v( \omega_k - C   )$ must be in $\Q$.
This terminates our evaluation of the sequence $(z_i)_{ i \ge 0}$ associated to $v$'s extension to $R[y ; \delta]$.
\item If $v(\omega_k) = v( C  ) = \frac{m_{k + 1}}{ n_{k + 1}  } \in\Q$ and $\overline{ x^{  \frac{m_{k + 1}}{ n_{k + 1}  }  }\omega_k     }   =  \overline{ x^{  \frac{m_{k + 1}}{ n_{k + 1}  }  }   C  }  $, 
$v( (  \Pi_{i = 1}^k x^{ \frac{m_i}{n_i}  }    ) z_{\ell}  ) > v( \omega_k  )  $. We continue with our evaluation by writing
\begin{align*}
\omega_{k + 1   } &= (  \Pi_{i = 1}^{k + 1} x^{ \frac{m_i}{n_i}  }    ) z_{\ell} ( \Pi_{i = 1}^{k + 1 } B_i   ) + C_3.
\end{align*}
\item  If $v(\omega_k) = v( C  ) = \frac{m_{k + 1}}{ n_{k + 1}  } \in\Q$ and $\overline{ x^{  \frac{m_{k + 1}}{ n_{k + 1}  }  }\omega_k     }   \neq  \overline{ x^{  \frac{m_{k + 1}}{ n_{k + 1}  }  }   C  }  $, then
$v( (  \Pi_{i = 1}^k x^{ \frac{m_i}{n_i}  }    ) z_{\ell}  )=  v( \omega_k  )  $
and $\gamma_{\ell + 1}  =  ( \overline{  x^{  \frac{m_{k + 1}}{ n_{k + 1}  }  }\omega_k}  -  \overline{x^{  \frac{m_{k + 1}}{ n_{k + 1}  }  }C     } )\Pi_{i = 1}^k \overline{B_i}^{-1} $. We write
\begin{align*}
\omega_{k  + 1 } &= (  \Pi_{i = 1}^{k + 1} x^{ \frac{m_i}{n_i}  }    ) z_{\ell + 1} ( \Pi_{i = 1}^{k + 1} B_i   ) + C_4.
\end{align*}
\end{enumerate}
For each $1 \le i \le 4$, $C_i$ is, as is $C$, an $\R$-linear sum of $D_{i,j}$ for various $i, j$. This is assured by Lemma \ref{construction three}. 

We point out that for each $\ell$, $v(z_{\ell})$ is determined in a finite number of steps. In case $v$ is determined by a finite sequence $( \omega_i  )^N_{i \ge - 1}$ for some $N \ge 0$, this is immediate. In the infinite case,
it follows from Lemma \ref{kompogoj} that $\lim_{k \rightarrow \infty} v( \omega_k ) = 0$, so, given that by Lemma \ref{construction two}, $v(C)$ is a sum of $v( \omega_i  )$, $v(\omega_k) < v(C)$ must hold for some $k$.
And since, in the step where $v( z_{\ell} )$ is determined, we get
$v( \Pi_{i = 1}^k  x^{ \frac{m_i}{n_i}  } z_{\ell}) = v( z_{\ell }  ) - \sum_{i = 1}^k v( \omega_{i - 1}   )   \le v( \omega_k )$ for each $\ell \ge 0$, we
get $v(z_{\ell }) \le \sum_{i = 0}^k v(  \omega_i) < 1$ by Lemma \ref{kompogoj}. Since for each $\ell \ge 0$, 
$z_{\ell + 1} = z_{\ell} - x^{- r_{\ell + 1}}\gamma_{\ell + 1}$, we did indeed find a unique sequence $(z_i)_{i \ge 0}$ that uniquely determines a valuation $v$ on $R[y ; \delta]$. This valuation is $v$'s extension from $\rweyl$ to $R[y ; \delta]$.
\endproof

The construction introduced in the proof of Theorem \ref{construction result} can be reversed. Given a valuation $v$ on $R[y ; \delta]$, we could use the reverse construction to find the sequence $( \omega_i)_{i \ge -1} \subseteq \rweyl$ associated to $v$'s restriction to $\rweyl$.

\section{Valuations on $\tilde{R}[y ; \delta  ]$}\label{grand}

The ring $\tilde{R}[y ; \delta  ]$ is an extension of $R[ y , \delta  ]  $ where $\tilde{R}$ is defined as $\R((x^{\Q}))$, the generalized power ring of sums $\sum_{q \in \Q} \alpha_q x^{-q}$ with well-ordered support. We  first show that every valuation on  ${R}[y ; \delta  ]$ can be easily extended to 
$\tilde{R}[y ; \delta  ]$.

\begin{lema} \label{rtildeext}
Every valuation on $R[y ; \delta]$ with residue field $\R$ can be extended to a valuation on $\tilde{R}[y ; \delta  ]$ with the same residue field. 
\end{lema}
\proof
Suppose first $v$ is defined on $R[y ; \delta]$ by a finite sequence $( z_i)_{i = 0}^k$ with $v(z_k) \not\in \Q$. Then, as we can write every 
$f \in \tilde{R}[y ; \delta]$ as $\sum_{i = 0}^n p_i(x) z_k^i$ with $p_i(x) \in \tilde{R}$, we define 
$v(f) = \min_{1 \le i \le n} \{v(p_i(x)) + i v(z_k)\}$. This gives us a well-defined valuation on $ \tilde{R}[y ; \delta]$ which clearly extends the one we defined on  $R[y ; \delta]$ in the previous section. 

Now suppose $v$ is defined by an infinite sequence $( z_i)_{i \ge 0}$ with $r = \lim_{i \rightarrow \infty} v(z_i) \le 1$. Define
 $z := y -  \sum_{i = 1}^{\infty} \alpha_i x^{- r_i} \in \tilde{R}[ y ; \delta]$. In the same way as before, write every $f \in \tilde{R}[y ; \delta]$ as $\sum_{i = 0}^n p_i(x) z^i$ with $p_i(x) \in \tilde{R}$ and define $v( f   ) = \min_{1 \le i \le n}\{ v(p_i (x))y + i (r - \mu  )  \}$ where $\mu$ is a positive infinitesimal. Thus we once more get $v$'s extension to 
$\tilde{R}[y ; \delta]$.
\endproof

In case $r = 1$, this is the only possible extension up to isomorphism of the value group, for $v(z)$ must be $1 - \mu$. This is because on the one hand, since
$v(z) > v(z_k) = r_{k + 1}$ for all $k \ge 0$, $v(z)$ is greater than any rational number $q < 1$. On the other hand, since $v[y - z, x] = v[y, x] = 0$ and the value group is commutative, $v(z) \le 1$. Thus if $v(z) \in \R$, $v(z) = 1$. But if we restricted $v$ to the quotient ring of the $\R$-algebra, generated by $x$ and  $z$, 
we would get a valuation on a division ring, isomorphic to $\mathcal{D}_1(\R)$, with a rational value group, residue field $\R$ and $v[ z, x] = v(zx  ),$ which contradicts Corollary \ref{ratvalgrp}.

\begin{trditev}\label{notrankone}
Let $v$ be a valuation on $\tilde{R}[y ; \delta  ]$ with residue field $\R$. Then the value group is not of rational rank one.
\end{trditev}
\proof
The only case in the proof of Lemma \ref{rtildeext} where it does not immediately follow that the value group is not $\Q$ is when $v$'s restriction to $R[y ; \delta]$ is constructed by an infinite sequence 
$( z_i )_{i  \ge 0}$ with $ r := \lim_{i \rightarrow \infty } v(y - z_i) < 1$. In this case, we can set $v(z) = r  \in \R$ and if $r \in \Q$, define $y^{(1)} :=  z$ and restart the construction of $v$. We may get another infinite sequence $(z_i^{(1)} )_{i \ge 0}$ with $r^{(1)} := \lim_{ i \rightarrow \infty } v(z_i^{(1)}) < 1$. If $r^{(1)} \in \Q$, we start over with $y^{(2)} :=  z^{(1)}$. Though we may have to repeat the process infinitely many times, the set 
$\{ z^{(j)}_k  \}_{j \ge 1, k \ge 0}$ is countable and $A := \{ v(  z^{(j)}_k)  \}_{j \ge 1, k \ge 0}$ is a well-ordered set of rational numbers smaller than one. At one point, $v(z)$ will have to be irrational for some $z = \sum_{q \in A} \alpha_q x^{- q} $ since we would otherwise get $z \in \tilde{R}$ such that $v( z) = 1$ which, as we have shown, contradicts the fact that the value group is rational.
\endproof

Recall that a pseudo-Cauchy sequence in a division ring $D$ with a valuation $v$ is a sequence $( a_{\lambda})_{\lambda \in \Lambda} \subseteq D $, where 
$\Lambda$ is an ordinal such that there exists $\lambda \in \Lambda$ for which 
$v( a_{\sigma} - a_{\rho} )  < v(a_{\rho} - a_{\tau}) $ for all  $\sigma, \rho, \tau \in \Lambda$  with $ \lambda \le \sigma < \rho < \tau$. 
Let $D'$ be an extension of $D$ and $v'$ and extension of $v$ to $D'$. Then $a \in D'$ is a limit of the pseudo-Cauchy sequence 
$( a_{\lambda})_{\lambda \in \Lambda}$ if $v'( a - a_{\sigma}) = v'(a_{\sigma + 1 } - a_{\sigma}  )$ for all $\sigma \in \Lambda$, $\lambda \le \sigma$.
As a byproduct of our investigations, we show that not every extension of a valued division ring by limits of pseudo-Cauchy sequences is immediate. This differs from the commutative case since, as Kaplansky proved in \cite{KAPL}, every extension of a valued field by limits of pseudo-Cauchy sequences is immediate.
\begin{posledica}
There exist division rings $D \subseteq D'$ and a valuation $v$ on $D$ which extends to a valuation $v'$ on $D'$ such that $D'$ is an extension of $D$ by limits of pseudo-Cauchy sequences in $D$ whereas $v'$ is not an immediate extension of $v$.
\end{posledica}
%
%
\proof
Let $v$ be a valuation on $R[y; \delta]$ with residue field $\R$ and value group $\Q$ as described in Section \ref{middle}. Then $v(x) = -1$ and $\tilde{R}[y ; \delta]$ is an extension of $R[y ; \delta]$ by limits of pseudo-Cauchy sequences. This holds because every
$ \sum_{i \in \Q} a_i x^{-i} \in \tilde{R} $ is a limit of the pseudo-Cauchy sequence $(  \sum_{i = 1}^k a_i x^{-q_i}  )_{k \ge 1}$ in $R[y ; \delta].$

As we have shown in this section, $v$ can be uniquely extended to $\tilde{R}[y ; \delta]$. By Proposition \ref{notrankone}, this extension is not immediate. Let $D$ be the quotient division ring of $R[y, \delta]$, to which $v$ uniquely extends, and ${D'}$ the quotient division ring of $\tilde{R}[y ; \delta]$, to which $v'$ uniquely extends, since both rings are Ore domains. $D'$ is not an immediate extension of the ring $D$ with valuation $v$, even though $D'$ is an extension of $D$ by limits of pseudo-Cauchy sequences.
\endproof

\section{Compatibility with orderings on $\rweyl$ and $R[y ; \delta]$} \label{orders}
In Section \ref{rweyl}, we mentioned that every strongly abelian valuation on a division ring with an ordered residue field is compatible with an ordering on the valued division ring. In this section, we will use a noncommutative version of the Baer-Krull theorem to determine all orderings on $\rweyl$ compatible with one of the valuations $v$ we have described in the previous sections. 
 We will then show which of these orderings on $\rweyl$ can be extended to an ordering on $R[y ; \delta]$ compatible with a $v$'s extension to $R[y ; \delta]$. 

Recall that an order $P$ on a division ring $D$ is compatible with a valuation $v$ on $D$ if for every $a, b\in D^*$ such that
$v(a) = v(b) < v(a - b)$, $ab \in P$ holds.

Let $v$ be a strongly abelian valuation on a division ring $D$ with a formally real residue field $\overline{D}$. 
Let $\Gamma$ be its value group. Let $s : \Gamma  \rightarrow D^* $ be a semisection of $v$, i.e., a map for which
\begin{enumerate}
\item $s(0) = 1$,
\item $v(s(g)) = g$ for all $g \in \Gamma$,
\item $s(g_1 + g_2) = s(g_1)s(g_2)u^2$ for some $u \in D^*$ for all $g_1, g_2 \in \Gamma$.
\end{enumerate}
Let $\chi : \Gamma / 2\Gamma \rightarrow  \{ -1, 1 \}$ be a group homomorphism called a character and let $\overline{P}$ be an ordering of $\overline{D}$.
Then, as it was shown in \cite{TSC},
\begin{align*}
P_{\chi} &= \{  a \in D  \mid \overline{a \cdot s(v(a))^{-1}} \chi( v(a) + 2\Gamma   ) \in \overline{P} \}
\end{align*}
is an order of $D$ compatible with $v$. Moreover, if $\overline{X}$ denotes all orders of the residue field, $X_v$ denotes all $v$-compatible orders on $D$ and $( \Gamma / 2\Gamma   )^*$ denotes the set of all characters on $\Gamma / 2\Gamma$, then by Proposition 3 of \cite{TSC}, the map
 \begin{align*} 
&f : \overline{X}  \times ( \Gamma / 2\Gamma   )^* \rightarrow X_v \\
&f ( (  \overline{P}, \chi)) = P_{\chi}
\end{align*}
is a bijection. The choice of a semisection $s$ on $\Gamma$ does not matter. Using $f$, we will now describe all orders on $\rweyl$ that are compatible with a valuation described in Section \ref{rweyl}.

Suppose $v$ is a valuation on $\rweyl$ associated to an infinite sequence $(\omega_i)_{i \ge - 1}$ with $\R$ as a residue field. There is only one possible order of $\R$, so the orders of $\rweyl$ compatible with $v$ will only depend on the characters 
$\chi : \Gamma / 2 \Gamma \rightarrow \{ -1, 1 \}.$ 
Then there are three different options for the value group $\Gamma \subseteq \Q \times \Z$.
\begin{enumerate}
\item $\Gamma$ is a 2-divisible subgroup of $\Q$,
\item $\Gamma$ is a non-2-divisible subgroup of $\Q$,
\item $\Gamma$ is a direct sum of a non-2-divisible subgroup of $\Q$ with $\Z$.
\end{enumerate}
Since in each case the value group $\Gamma$ is generated by $ \{ v(\omega_i) \}_{i \ge -1}$, the characters and thus the $v$-compatible orders will be determined by the signs of the  $\omega_i$.

In the first two cases, $v$ is determined by an infinite sequence $(\omega_i)_{i \ge -1}$ with $v( \omega_i ) = \frac{m_{i + 1}}{  n_{i + 1}  }   \in \Q$ for each 
$i \ge - 1$. In the third case, $v$ is determined by a finite sequence  $(\omega_i)_{i = -1}^N$ with $v( \omega_i ) = \frac{m_{i + 1}}{  n_{i + 1}  }   \in \Q$ for each $-1 \le i \le N  - 1$ and $v( \omega_N) \not\in \Q$.

If the value group is a $2$-divisible subgroup of $\Q$, then for each $\omega_i$, there is a $j \neq i$ such that $K_{i,j} \in \Z$ is odd while $K_{j, i} \in \Z$ is even. Since 
$\overline{ \omega_{i-1}^{K_{i,j}} \omega_{j-1}^{K_{j,i}    }  } =  \alpha_{i,j}  $, it follows that $\omega_{i -1} > 0$ if and only if 
$ \alpha_{i,j} > 0$. 

Conversely, if the value group is a non-$2$-divisible subgroup of $\Q$, we choose one $i$ such that $ n_{i } $ is divisible by the greatest power of two that divides $n_{j } $ for any $j \ge 0$. After choosing either $\omega_i < 0$ or $\omega_i > 0$, the order on $\rweyl$ is defined. 

In the last case, where the value group is a direct sum of a non-2-divisible subgroup of $\Q$ and $\Z$, the order is determined by
choosing either $\omega_i < 0$ or $\omega_i > 0$ and, independently, either $\omega_N > 0$ or $\omega_N < 0$ where $i$ is as in the second case and $v(\omega_N ) \not\in \Q$. All four combinations define an ordering on $\rweyl$.

We have thus proven the following proposition. 
\begin{trditev}
Suppose $v$ is a valuation on $\rweyl$ with residue field $\R$ and value group $\Gamma$. Then:
\begin{enumerate}
\item If $\Gamma$ is a $2$-divisible subgroup of $\Q$, there is a unique $v$-compatible ordering on $\rweyl$.
\item If $\Gamma$ is a non-$2$-divisible subgroup of $\Q$, there are two $v$-compatible orderings on $\rweyl$.
\item If $\Gamma$ is a direct sum of a non-$2$-divisible subgroup of $\Q$ with $\Z$, there are four possible $v$-compatible orderings on $\rweyl$.
\end{enumerate}
\end{trditev}
\subsection{Extensions of orderings on $\rweyl$ to orders on $R[y ; \delta]$} 

In this section, we show which orders on $\rweyl$ are extendable to an order on $R[y ; \delta]$, thereby answering the question posed by Marshall and Zhang in \cite{MA1}. 

Every order on $\rweyl$ is compatible with a unique finest valuation $v$ on the same ring with residue field $\R$, as proved in \cite{MA1}.
Suppose $v$ is a valuation on $\rweyl $ associated to an infinite sequence $( \omega_i  )_{i \ge - 1}$ with 
$v( \omega_{i-1}) = \frac{ m_{i}  }{ n_{i } }$ and 
$\omega_{i } = x^{m_{i }} \omega_{i-1}^{n_{i}} - \beta_{i}$.
In Section \ref{middle}, we showed that provided both conditions of Theorem \ref{construction result} are fulfilled, $v$ can be uniquely extended to a valuation $v'$ on $R[y; \delta ]$ with residue field $\R$ if the value group is $2$-divisible and that it has two extensions to $R[y ; \delta]$ with the same residue field if the value group is non-$2$-divisible.
 Here we show all $v$-compatible orders on $\rweyl$ we have described in the first part of this
section can be extended to a $v'$-compatible order on $R[y ; \delta]$ for some extension $v'$ of $v$ from $\rweyl$ to $R[y ; \delta]$.

\begin{izrek}
Let $P$ be an ordering on $\rweyl$ and $v$ be the unique finest valuation on $\rweyl$ compatible with $P$. Then: 
\begin{enumerate}
\item The order $P$ can be extended to an ordering on $R[y ; \delta]$ if and only if $v$ can be extended to a valuation on $R[y ; \delta]$ with residue field $\R$. 
\item If the $v$-value group $\Gamma$ is a $2$-divisible subgroup of $\Q$, then the extension of $P$ is unique. If on the other hand, 
$\Gamma$ is a subgroup of $\Gamma_1 \times \Z$, where $\Gamma_1$ is a non-$2$-divisible subgroup of
$\Q$, i.e., when $\Gamma$ is not a $2$-divisible subgroup of $\Q$, there are two extensions of $P$ to $R[y ; \delta]$. Each of the two extensions of $v$ to a valuation $v'$ on $R[y ; \delta]$ with residue field $\R$ uniquely determines one of the two of $P$'s extensions to $R[y ; \delta]$. 
\end{enumerate}
\end{izrek}
\proof
The first statement of the theorem follows from the fact that every ordering on $R[y ; \delta]$ is compatible with a valuation on the same ring with residue field $\R$. 

To prove the second statement, suppose that $v$ is the unique $P$-compatible valuation on $\rweyl$ with residue field $\R$ that extends to a valuation on 
$R[y ; \delta]$ with the same residue field.

If the value group $\Gamma$ of $v$ on $\rweyl$ is either a $2$-divisible or non-$2$-divisible subgroup of $\Q$, then the value group $\Gamma'$ of $v$'s extension to $R[y ; \delta]$ is $ \Q$. 
In this case, there is exactly one $v'$-compatible order of $R[ y ; \delta ]$ for each of $v$'s extensions $v'$ to $R[y ; \delta]$. 

Suppose $\Gamma$ is a $2$-divisible subgroup of $\Q$. Then there is a unique extension $v'$ of $v$ to $R[y ; \delta]$. 
It follows that in case $\Gamma$ is a $2$-divisible subgroup of $\Q$, the only $v$-compatible order on $\rweyl$ extends to an order of 
$R[y ; \delta]$ that is compatible with $v'$. 

If, on the other hand, $\Gamma$ is a non-2-divisible  subgroup of $\Q$, there are two extensions $v'$ of $v$ to $R[y ; \delta]$. 
We will now show that for each of the 
$v$-compatible orderings on $\rweyl$, there is a unique extension $v'$ of $v$ to $R[y ; \delta]$ such that the ordering on $\rweyl$ can be extended to the unique $v'$-compatible ordering on $R[y ; \delta]$. 

In this case, a $v$-compatible ordering on $\rweyl$ is, as we have shown in the beginning of this section, uniquely determined by the sign of $\omega_{i - 1}$ where $i \ge 1$ is such that $n_i$ is divisible by the greatest power of two that divides $n_j$ for any $ j \ge 1$. Furthermore, $v'$, the extension of $v$ to $R[y ; \delta]$, is uniquely determined by choosing the sign of $\tilde{\gamma_i}$ for this $i$. 

We first choose an extension $v'$ of $v$ to $R[y ; \delta]$. We observe that $x^{ \frac{m}{n} } > 0$ must hold for every $\frac{m}{n} \in \Q$ since all rational powers of $x$ are in $R[y ; \delta ]$. 
Since $\overline{ x^{ \frac{ m_{i}}{n_i}   } \omega_{i - 1}  } = \tilde{\gamma_{i }} $ for each $i \ge 1$, 
$\tilde{\gamma_{i }}\omega_{i -1 }  > 0$ must hold for all $i \ge 0$ for the order to be extendable to a $v$-compatible order on $R[y ; \delta]$. 
This holds for exactly one of the two $v$-compatible orders on $\rweyl$.  It is clear from the construction that for each ordering on $\rweyl$, there is exactly one extension $v'$ of $v$ to $R[y ; \delta]$ such that this ordering is extendable to the unique $v'$-compatible ordering on $R[y ; \delta]$.


In case $\Gamma$ is a subgroup of $\Q \times \Z$ of rational rank two, $\Gamma' = \Q \times \Z$ holds. In this case, there are two $v'$-compatible orderings on $R[y ; \delta]$ for every extension $v'$ of $v$ to $R[y ; \delta]$. We will now show that for each of the four $v$-compatible orders on $\rweyl$, there is a unique extension $v'$ of $v$ to $R[y ; \delta]$ and a unique $v'$-compatible ordering $P'$ on $R[y ; \delta]$ such that $P'$ is an extension of $P$.
The ordering $P$ compatible to a valuation $v$ on $\rweyl$ is determined by the signs of $\omega_{i-1}$ and $\omega_N$ where $i \ge 1$ is such that $n_i$ is divisible by the greatest power of two that divides $n_j$ for any $ j \ge 1$, and $v(\omega_N) \not\in \Q$.
The extension $v'$ of $v$ to $R[y ; \delta]$ and the $v'$-compatible ordering on $R[y ; \delta]$ that extends $P$ are the valuation $v'$ for which $\omega_{i - 1}\tilde{\gamma_i} > 0$ and the $v'$-compatible ordering that agrees with the signs of $\omega_{i-1}$ and $\omega_N$ in $P$.
%

We have thus proved the second statement of the theorem.
\endproof

\end{document}